\documentclass{article}

\usepackage[utf8x]{inputenc}
\usepackage[T1]{fontenc}
\usepackage{lmodern, textcomp}
\usepackage[english]{babel}

\usepackage{fancyhdr}
\usepackage[pagewise]{lineno}
\usepackage{import}

\usepackage{amssymb}
\usepackage{amsmath}

\usepackage{graphicx}
\graphicspath{{figs/}}
\usepackage{pgfplots, pgfplotstable, booktabs}
\pgfplotsset{compat=1.8}

\newcommand{\eps}{\varepsilon}
\newcommand{\R}{\mathbb{R}}
\newcommand{\NN}{\text{NN}}

\title{A neural network closure for the Euler-Poisson system based on kinetic simulations}
\author{L. Bois$^{1,2}$, E. Franck$^{1,2}$, L. Navoret$^{1,2}$, V. Vigon$^{1}$}
\date{August 2020}

\begin{document}

\maketitle

\centerline{1. Institut de Recherche Math\'ematique Avanc\'ee, UMR 7501,}
\centerline{ Universit\'e de Strasbourg et CNRS, 7 rue Ren\'e Descartes,}
\centerline{ 67000 Strasbourg, France}
\bigskip

\centerline{2. INRIA Nancy-Grand Est, TONUS Project, Strasbourg, France}

\setcounter{tocdepth}{2}
\begin{abstract} This work deals with the modeling of plasmas, which are charged-particle fluids. Thanks to machine leaning, we construct a closure for the one-dimensional Euler-Poisson system valid for a wide range of collision regimes. This closure, based on a fully convolutional neural network called V-net, takes as input the whole spatial density, mean velocity and temperature and predicts as output the whole heat flux. It is learned from data coming from kinetic simulations of the Vlasov-Poisson equations. Data generation and preprocessings are designed to ensure an almost uniform accuracy over the chosen range of Knudsen numbers (which parametrize collision regimes). Finally, several numerical tests are carried out to assess validity and flexibility of the whole pipeline.
\end{abstract}
\vspace{0.2cm}

\vspace{0.2cm}


\tableofcontents

\section{Introduction}

Plasmas are gases composed of charged particles, i.e. ions and electrons, which are the subject of many studies because they are a very common state of matter in the universe (e.g. in stars, ionosphere). They are also present in industrial devices, like in fusion reactors. Plasma dynamics is complex since particles have both long range interactions, through the self-consistent electromagnetic fields, and short-range collisions. The  most complete model to describe plasmas is given by the Vlasov equation, which is a kinetic equation satisfied by the distribution function of particles in the position-velocity $(x,v)$-phase space and which is coupled with the electromagnetic equations. To reduce the dimention, it is possible to use \emph{fluid models} which are derived from the kinetic Vlasov equation with a collision operator by taking the velocity moments corresponding to the monomials $1$, $v$, $v^2$. This raises up three equations linking the four following quantities: the \emph{density} (the zeroth order moment), the \emph{mean velocity} (involving the two first moments), the \emph{temperature}  (involving the three first order moments) and  the \emph{heat flow}  (involving the four first moments).
Fluid models are very attractive as they are much computationally cheaper as they involve only spatial quantities. However, they require an additional equation called a \emph{closure} to link the heat flow to the other three moments. 
 
 The first two fluid models obtained in the collisional regime are the Euler and the Navier-Stokes system. The closure is obtained from an asymptotic Chapman-Enskog procedure, assuming that the distribution function is close to be at thermodynamical equilibrium due to collisions. These equations are thus valid for small Knudsen number $\varepsilon > 0$, which is the mean free path between two collisions divided by the characteristic length studied. The Navier-Stokes system provides $O(\varepsilon)$ corrections of the pressure tensor and heat flux. Higher order corrections have been considered but lead to ill-posed systems (Burnett equations). We refer to \cite{Braginskii_1963} for a specific derivation in the case of plasmas and to \cite{degond} for a review.
 
To extend the validity of fluid models for larger Knudsen numbers, larger moment systems have been proposed for the Boltzmann kinetic equation, with neutral particles: the Grad 13 order moment model using perturbative theory \cite{Grad}, the Levermore 14 order moment model using entropy maximisation closure \cite{Levermore, Levermore1998}. They both suffer from intrinsic mathematical imperfections: a lack of hyperbolicity for the former \cite{CaiFanLi}, so-called realizibility issues for the latter \cite{Junk, Schneider}. However, recent developments propose corrections for the Grad model \cite{CaiFanLi} and the entropy methodology and has been very successful in the context of radiative transfer with compact velocity domain \cite{DubrocaFeugeas, GarrettHauck}. 

Specific closures have been developed for plasma fluid models. In electrostatic regime, the Hammett-Perkins closure \cite{Hammett, Hammett1992} is designed to recover the Landau damping effect, which results from the transfer from spatial modes to velocity modes and which leads to the damping of the electrostatic energy. To this end, the heat flux depends on the temperature through a non-local integral dissipative operator. Dissipation is also introduced in the momentum equation. Let us mention the recent work \cite{Manfredi} for a discussion on such closure. When considering magnetized plasma, several closures have also been proposed for fluid models, like for Magnetohydrodynamics (MHD) or gyrofluid equations. See for instance \cite{1957_Chew, Snyder, Brizard_gyrofluid, Negulescu2016}. Recently, fluid closures that preserve the Hamiltonian structure of the Vlasov equation have also been developed \cite{Perin_Hamiltonian2015, Tassi}. 

Other strategies consist in adding kinetic effects. Let us mention model reduction technic to obtain fluid models like the water-bag method \cite{Besse2009, BesseGyro, Desjardin_momentmethod}, where the kinetic distribution function is approximated with piece-wise constant functions in velocity. Kinetic effects can also be numerically introduced through micro-macro decomposition methods \cite{CrouseillesDegond, DegondDimarcoMieussens, Crestetto}.

A more recent approach consists in using supervised machine learning to find a data-driven closure. Artificial neural networks (ANN) have proven very efficient to interpolate data and detect underlying structures. In particular, convolutional neural networks are very efficient to analyse images and thus the outputs of numerical simulations. There are numerous works for designing physics-based models using neural networks (see for instance \cite{BeckFladMunz, WangWuXiao_Reynolds, Duraisamy2019}). Such tools have already been used in two works concerning moment closure. In \cite{han}, in the context of neutral gases, the authors propose to learn the appropriate higher moments required in the model. In \cite{maulik}, which concerns plasmas models, the closure is directly learned from the analytic closure like the Hamett-Perkins one. 

\paragraph{Description of our work.}

Here, we introduce a data-driven closure, based on a fully-convolutionnal neural network, which is valid for a large spectrum of collisional regimes (i.e. a large range of $\varepsilon$). 
Data are obtained from numerical simulations of the Vlasov-Poisson system satisfied by the distribution function and the electric potential. This method has already been mentioned in \cite{han} for neutral gases, but the authors preferred to develop another approach. Here we further investigate the design of the closure to obtain accurate predictions.

Collisions between charged particles would normally be modeled with the Fokker-Planck-Landau operator. However, as it is usually done, we replace it by a BGK (Bhatnagar-Gross-Krook) operator, that models the relaxation of the distribution function towards the equilibrium distributions, the Maxwellians. The inverse of the Knudsen number $\varepsilon$ is in prefactor of this operator. 

Knudsen numbers $\varepsilon$ range is chosen equal to $[0.01,1]$. Indeed, for $\varepsilon < 0.01$, the Navier-Stokes closure already gives good results. 
An upper bound of the Knudsen interval has to be prescribed. Here we choose $\varepsilon = 1$ for embrassing mildly collisional regimes, but larger values could be considered. 

The closure takes as input three one-dimensional vectors corresponding to the spatial discretized density, mean velocity and temperature, and one vector corresponding to the parameter $\varepsilon$. The  closure returns an estimation of the spatial discretized heat flux.

The core of the closure is a fully convolutional neural network which has a \emph{V-Net} architecture. It was first described and developed by Milletari et al. \cite{milletari} to treat medical images (3D signals). It is an evolution of the U-Net, developed for biomedical 2D images \cite{ronneberger}. It consists in a succession of several convolution kernels, down-samplings and up-samplings that perform multi-scale analysis of the signal. No fully connected layer is used, so the whole network has relatively few parameters. 
The architecture of the network can mostly be described by three hyperparameters: the number of levels of the "V", the depth (which rules the number of channels at each level), and the size of the kernels of convolution.  The influence of these parameters on the performance is investigated.

As usual in neural networks construction, some processing of the data is required to predict the heat flux in the widest possible range. Therefore multiple steps of processing are included in the closure, on top of the central neural network.
One of these is the normalization of the heat flux with the Navier-Stokes approximation, that helps lowering the otherwise large relative error on predictions of low heat fluxes.
Another weakens the dependency of the closure on the underlying mesh size, by performing a data slicing before applying the neural network part.
We also use a resampling strategy for discretizations with a resolution different from the one used to train the neural network.

The insertion of this closure in a fluid solver raises mathematical issues. Indeed, it is not guaranteed that the neural network closure is dissipative and thus it could lead to instabilities. To prevent such scenario, a smoothing of the prediction is added into the closure. 

The paper is organised as follows. In Section \ref{sec:fluidclosure}, we introduce the Vlasov-Poisson system and the moment closure issue. In particular, we give the Euler and Navier-Stokes closures. Then the neural network closure strategy is explained in Section \ref{sec:networks}: prediction strategy, architecture of the network and training methodology are presented. Then Section \ref{sec:data} details the data generation and Section \ref{sec:processing} their processing. Finally, in Section \ref{sec:results}, we carry out several numerical tests to quantify and analyze the accuracy of the closure.


\section{Fluid closure for Vlasov-Poisson}
\label{sec:fluidclosure}

In this section we present the kinetic description of the plasma dynamics in one space dimension and its fluid approximations.

\subsection{Kinetic model}
\label{sec:kinetic}

A kinetic model is a model interested in the function $f:(x,v,t) \mapsto f(x,v,t)$ that describes the evolution of the distribution of the particles in the $(x,v)$-phase space, where $x\in [0,L]$ denotes the space variable, $v\in\mathbb{R}$ the velocity variable and $t \in \R$ the time. We will consider periodic boundary conditions in space. From this distribution function can be computed some macroscopic physical quantities. The first three moments give the particle density $\rho(x,t)$, the mean velocity $u(x,t)$ and the total energy $w(x,t)$ defined as:
\begin{align}
&\rho(x,t) = \int_{\R} f(x,v,t) dv,\quad \rho(x,t) u(x,t) = \int_{\R} f(x,v,t) v dv,\label{def:rho-u}\\
&w (x,t) = \frac{1}{2} \int_{\R} f(x,v,t) v^2 dv.\label{def:w}
\end{align}
We can also define the pressure $p(x,t)$, the temperature $T(x,t)$ and the heat flux $q(x,t)$:
\begin{align}
&p(x,t) = \int_{\R} f(x,v,t)(v-u(x,t))^2 dv,\quad \rho(x,t) T(x,t) = p(x,t),\label{def:p-T}\\
&q(x,t) = \int_{\R} \frac{1}{2} f(x,v,t) (v-u(x,t))^3 dv.\label{def:q}
\end{align}
Note that we have the following relation:
$w = \rho u^2 /2 +  \rho T /2 = \rho u^2 /2 +  p /2$.\\
The evolution of the distribution is described by the Vlasov equation:
\begin{align}
    \label{eq:vlasov}
    \partial_t f
    +\ v \partial_x f
    -\ E \partial_v f
    =\ Q(f),
\end{align}
where $E(x,t)$ is the self-induced electric field, which satisfies the Poisson equation:
\begin{align}
    \label{eq:poisson}
    E = - \partial_x \phi \quad,\quad \partial_{xx} \phi = \rho - \int_{[0,L]} \rho\,dx.
\end{align}
Here $\phi(x,t)$ denotes the electric potential.

The source term $Q$ is called a collision operator and allows the model to take into account the collisions between particles. Different collision operators can be considered to deal with different situations. In this work we use the BGK operator (Bhatnagar, Gross and Krook), built to conserve the mass, momentum and kinetic energy of the system, and model the relaxation induced by the collisions toward a local equilibrium distribution $M(f)$.
This operator simply reads
\[ Q(f) = \frac{1}{\varepsilon} ( M(f) - f ), \]
where $M(f)$ is called the Maxwellian of $f$ and is given by
\begin{equation}
    \label{eq:maxwellian}
    M(f)(x, v, t)
    = \frac{\rho(x,t)}{\sqrt{2\pi T(x,t)}}\ e^{-\frac{(v-u(x,t))^2}{2 T(x,t)}},
\end{equation}
with $\rho$, $u$ and $T$ the density, mean velocity and temperature associated to $f$ and defined in \eqref{def:rho-u}-\eqref{def:p-T}.
The parameter $\varepsilon > 0$ is called the Knudsen number and represents the mean free path between two collisions divided by the characteristic length studied. In the limit $\eps \to 0$, the distribution function is expected to be closed to local equilibria. In this regime, the phase-space dynamics could be infered from the spatial dynamics of its macroscopic moments $\rho$, $u$, $T$.

The kinetic equation \eqref{eq:vlasov} can be easily solved numerically in this one dimensional case. Several numerical methods have been proposed in the literature \cite{sonnendrucker,dimarco2014}. In this work, we use a Finite Difference/Finite Volume method similar to the one introduced in \cite{PhamHelluy, helluy} and described in appendix \ref{sec:vlasov-num}. Although it is possible to extend such computations in dimension 2 or 3, the computational cost becomes very prohibitive. We are therefore led to consider fluid models.

\subsection{Fluid model}
\label{sec:fluid}

A fluid model of a plasma describes the evolution of its density $\rho(x,t)$, mean velocity $u(x,t)$ and kinetic energy $w(x,t)$, instead of the distribution of its particles $f(x,v,t)$ in a kinetic model. Such models are far less expensive to simulate numerically than the original kinetic model. Indeed, this requires the computation of only three spatial quantities ($3\times N_x$ unknowns), instead of the full kinetic distribution ($N_x\times N_v$ unknowns), where $N_x$ and $N_v$ stands for the number of discretization points in space and velocity.

Fluid models can be obtained by using the moment method \cite{degond}. Formally, it consists in computing the first three moments in velocity of the Vlasov equation, i.e. multiplying Equation \eqref{eq:vlasov} by $1$, $v$ and $v^2/2$ and then integrating in velocity:
\begin{align*}
\text{for $p = 0$, $1$, $2$,}\quad \int_{\R} v^p
    \left(\partial_t f
    +\ v \partial_x f
    -\ E \partial_v f\right)\, dv
    =\ \int_{\R} v^p Q(f)\, dv.
\end{align*}
Since the collision operator conserves mass, momentum and energy, the right-hand sides vanish. Therefore, it results in the following system:
\begin{equation}
    \label{eq:fluid}
    \left\{
    \begin{aligned}
        \partial_t \rho +  \partial_x (\rho u) &= 0, \\
        \partial_t (\rho u) + \partial_x (\rho u^2 + p) &= -E \rho, \\
        \partial_t w + \partial_x (w u + pu + q) &= - E \rho u,\\
    \end{aligned}
    \right.
\end{equation}
where $p$ is the pressure and $q$ the heat flux defined in \eqref{def:p-T}-\eqref{def:q}. Note that the pressure $p$ is actually a function of $\rho,  u, w$ since we have the following relation: 
$$p = 2 w -  \rho u^2.$$ The electric field $E$ is still given by the Poisson equation \eqref{eq:poisson}.

We thus get a system of three equations on the four variables $\rho$, $u$, $w$ and $q$.
For it to be closed, we need a fourth equation connecting these four unknowns, called a closure. Usually this closure consists in replacing the true heat flux $q$ with a simplified one $\hat q$  given as a function of the other quantities , and that can be written
\[
\hat q = \mathcal{C}(\varepsilon, \rho, u, T).
\]
Note that this closure depend on the physical regime we consider through the parameter $\eps$. Let us also mention here that in higher dimension, the pressure stress tensor is not completely determined by $\rho$, $u$ and $w$ and an additionnal closure relation is necessary.

In fluid regimes, where most of the kinetic effects can be neglected, two closures are classically considered: the Euler or the Navier-Stokes closures. The Euler closure is valid in regimes where the distribution function is closed to be Maxwellian: $f = M(f) + O(\varepsilon)$. Using this ansatz, we obtain the Euler closure: \[
\hat q  = 0.
\]  When we are interested in regime with $O(\eps)$ deviation from the Maxwellians, $f = M(f) + \varepsilon g + O(\varepsilon^2 )$, we get the Navier-Stokes closure: \begin{equation}
\label{eq:ns}
\hat q  = -\frac{3}{2} \varepsilon\, p\, \partial_x T.
\end{equation} Note that it is a $O(\eps)$ correction of the Euler closure. We refer to \cite{degond} for more details.

More complex closures have been developped to capture more kinetic effects, in regimes where the disbribution function is more distant from the Maxwellians.
Such closures are chosen non local, meaning that the heat flux $q(x)$ at location $x$ does not rely on the values of $\rho$, $u$, $T$ and their derivatives at location $x$ only, but on their values on all the domain. Indeed, the true heat flux, given by a third order moment of the kinetic distribution (Eq. \eqref{def:q}), depends on the full distribution in velocity at location $x$ and thus retains non-local information about the dynamics.

The goal of our work is to provide a method to provide a new closure, where the function $\mathcal{C}$ is implemented using a neural network. The next section describes this approach in more details.

\section{Closure with a neural network}
\label{sec:networks}

In this section we introduce the overall principle of our method to build the fluid closure presented in Section \ref{sec:fluid} with a neural network. First we describe the basic functioning of the neural network, before introducing the other operations included in the closure, and that surround the neural network. We end this section with the detailed description of our neural network and of its training.

\subsection{Interpolation of the heat flux with a neural network}

In this work we are interested in the ability of machine learning to find a function $C$ that maps the Knudsen number, the density, the mean velocity and the temperature of the plasma to its heat flux. Such a mapping cannot be exact as the heat flux is not a function of these four quantities but a function of the particles' distribution in the phase space. Nonetheless, the ability of neural networks to find non obvious patterns and correlations can be used to provide a good approximation of the heat flux, that can then be used as closure for the fluid model described in Section \ref{sec:fluid}.

Neural networks enable to build functions relying on many parameters. The obtained closure can be formally written as:
\[
\hat q = C_{\hat\theta}(\varepsilon, \rho, u, T),
\]
where $\hat\theta$ denotes the set of parameters of the obtained network. Here the closure is chosen to be \textbf{non-local}: $\rho$, $u$, $T$ as well as $q$ are all spatial discretized quantities.
The closure takes as input four vectors or, in the neural network terminology, a 1D signal with four channels corresponding respectively to the Knudsen number (turned into a constant vector), the density, the mean velocity and the temperature at spatial discretization points:
$$X = (\eps, \rho, u, T) \in (\R^{N_x})^4,$$
where $N_x$ is the number of spatial discretization points. The output of the closure is one vector of size $N_x$ 
$$ Y = C_\theta (\varepsilon, \rho, u, T) \in \R^{N_x},
$$
that will correspond to the estimated heat flux $\hat q \in \R^{N_x}$.

To find such a function, there are two main steps:
\begin{enumerate}
    \item First, by setting the \emph{architecture} of the neural network, we define a family of functions $C_\theta$ parametrized by $\theta \in \Theta$, where $\Theta$ represents the set of all the possible parameters of the neural network. 
    \item  Then, by \emph{training} the neural network, we find a set of parameters $\hat \theta$ so as to minimize the error of the neural network predictions on a dataset, for which the true heat flux is known.
\end{enumerate}

The neural network is the core component of the closure $C_\theta$, but it is not its only one. Indeed, for the data to be usable by the neural network and to provide satisfying results, it needs to go through a certain amount of processing. Section \ref{sec:integration} introduces these other components and their articulation in the closure. Sections \ref{sec:architecture} and \ref{sec:training} then describe respectively the architecture and the training of the neural network.

\subsection{Detailed composition of the closure}
\label{sec:integration}

In practice, the closure is used at each iteration in time to compute the heat flux $q$ on all the mesh from the input $(\varepsilon, \rho, u, T)$ on all the mesh.
Such a use of a neural network based closure into a numerical scheme requires some work to transform the data handled by the solver into data usable by the network, and vice-versa.
This section briefly introduces the different steps involved in the closure and their role in relation to the network.

\begin{figure}
    \centering
		\includegraphics{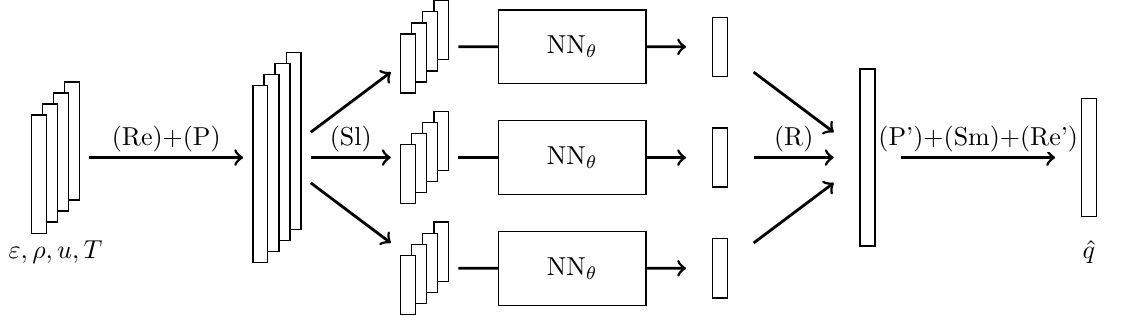}
    \caption{Graphic of the composition of the closure. The different operations are (Re)-(Re') resampling, (P) pre-processing, (Sl) slicing, ($\NN_\theta$) neural network, (R) reconstruction, (P') post-processing and (Sm) smoothing.}
    \label{fig:nn-integration}
\end{figure}

There are three main processing operations designed to integrate the closure into the numerical scheme. The first one deals with the difference in resolution between the numerical scheme and the data used to train the network, by resampling the input of the network to its training resolution, and its output back to the original resolution. The second one deals with the difference in length between the resampled data of the scheme and the inputs handled by the network, by slicing the data into several pieces to give to the network, and then aggregating the outputs to reconstruct the full heat flux. The third one deals with the stability of the resulting numerical scheme, by smoothing the data from the network to prevent any oscillation to propagate into the numerical scheme. These three steps are part of a whole process, that can be broken downs as follows:
\begin{equation*}
  C_\theta:\quad X \enskip \overset{\text{(Re)}+\text{(P)}}{\longmapsto} \enskip \tilde X\enskip \overset{\text{(Sl)}}{\longmapsto} \enskip (\tilde X_j)_j\enskip \overset{(\NN_\theta)}{\longmapsto}\enskip (\tilde Y_j)_j\enskip \overset{\text{(R)}}{\longmapsto} \enskip \tilde Y\enskip \overset{\text{(P')}+\text{(Sm)}+\text{(Re)}}{\longmapsto} Y.
\end{equation*}
This process is illustrated by figure \ref{fig:nn-integration}, and the different steps are described below.

\begin{description}
    \item[1. Resampling (Re) and pre-processing (P)] The input $X$ of size $N_x$ is resampled with a Fourier method (that relies on the periodicity of the signal) and pre-processed into a signal $\tilde X$ of size $N_x'$. As explained in section \ref{sec:res-resolution}, the resampling changes the resolution of the input to match the one used for training and get better results. The pre-processing step that follows is a common step in the use of neural networks that allows to improve the accuracy. It is detailed in section \ref{sec:in-std}.
    
    \item[2. Slicing (Sl)] The inputs are sliced into several pieces $(\tilde X_j)_j$ of size $N$ before being handed over to the network, where $N$ is the size of the input the neural network works with.
    
    \item[3. Neural network (NN$_{\theta}$)] For each piece of input $\tilde X_j$, the network predicts a piece of heat flux $\tilde Y_j$, of size $N$.
    
    \item[4. Reconstructing (R)] The signals $\tilde Y_j$ are aggregated in order to reconstruct the whole heat flux. An added benefit of this slicing and reconstructing process is that it allows to prevent some edge effect introduced by the network: by using overlapping pieces, the ends of each piece can be ignored when reconstructing the whole signal. As described in section \ref{sec:slicing}, the amount of overlapping is a parameter that can be tweaked.
    
    \item[5. Post-processing  (P'), smoothing (Sm) and resampling (Re')] \hfill \\
    The post-processing is an operation of normalization designed to improve the accuracy of the network and detailed in section \ref{sec:out-norm}. The smoothing of the signal is here to avoid oscillations coming from the predictions to propagate and amplify, thus making the numerical scheme unstable. The intensity of this smoothing can be tweaked, and is discussed in section \ref{sec:res-smoothing}. Finally it is followed by a resampling to recover the original resolution, for the output to be used in the rest of the computations.
\end{description}
These different operations will be fully exposed in Section \ref{sec:processing}.
We now turn to the heart of the closure: the neural network.

\subsection{Architecture of the neural network}
\label{sec:architecture}

The closure $C_\theta$ depends on the architecture of the network.
In our work we use a 1D version of the V-Net \cite{milletari, ronneberger} implemented with the Tensorflow 2 library through its python interface.

\begin{figure}
    \centering
    \includegraphics{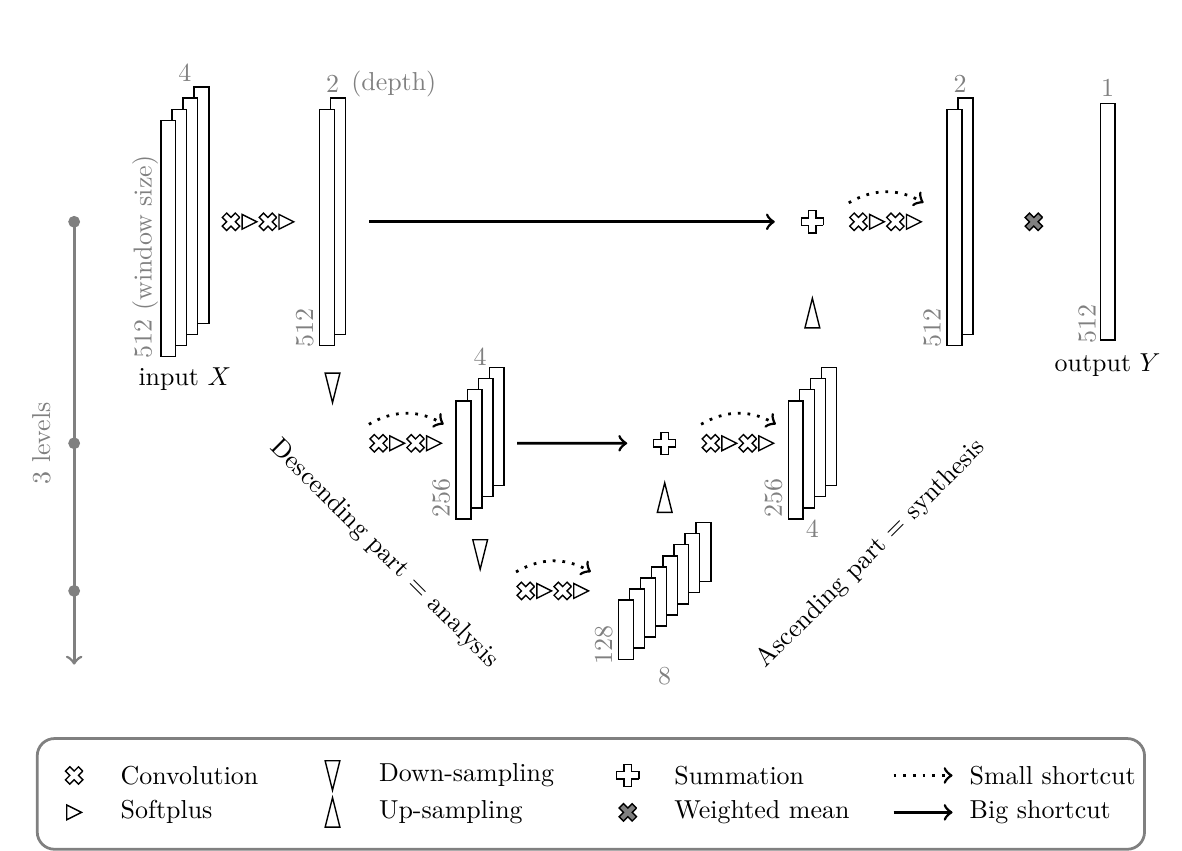}
    \caption{Graphic of a 1D V-Net with a window of size $N=512$, a depth $k = 2$ and $\ell = 3$ levels.}
    \label{fig:v-net}
\end{figure}

The V-net is a fully convolutional neural network meaning that it is based only on a successive sequence of convolutions (no fully connected layers). It consists of two parts. First a descending part that acts as an "encoder", and decomposes the input $X$ into multiple features, performing a multiscale analysis. Then an ascending part that acts as a "decoder", and synthesizes the created features to predict the output $Y$. As we will see below, it can be characterized by three hyperparameters: the number of levels $\ell$, the depth $d$ of the first convolution and the kernel size $p$ of all convolutions. The composition of the V-net illustrated Figure \ref{fig:v-net} is the following:
\begin{description}
\item [Input] As described previously, it is a 1D signal $X$ made of 4 channels $(\varepsilon, \rho, u, T)$, and with a given length called the \emph{window size} that we set to $N$.

\item [Initialization] We first perform successively two 1D convolutions%
\footnote{For the sake of completeness, we recall the formula for a one dimensional convolution: applied on an input $X$ of shape $(N,d)$ to get an output $Y$ of shape $(N,d')$, a 1D convolution with an odd kernel size $p$ uses a kernel $K$ of shape $(p,d,d')$ and we have
$$
Y_{i,k} = \sum_{j=1}^d \sum_{di=1}^p \tilde{X}_{i+di, j} K_{di,j,k},
$$
where $\tilde{X}$ is the input $X$ padded on both ends so that $Y{i,k}$ is well defined for $i$ up to $N$.
}
with a kernel size of $p$ (we tried $p=5, 7, 9, 11$), both followed by the softplus activation function $s(x) = \ln(1+e^x)$, and change the number of channels from 4 to a depth $d$ (we tried $d=4,5,6,8$).  All convolutions preserve the size of the signal with the help of a constant padding that extends the signal by continuity. The original V-Net uses the ReLu activation function, but the regularity of the softplus activation function seems to give us better results.

\item [Descending part] The V-net is made of $\ell$ levels (we tried $\ell=3, 4, 5$), so $\ell-1$ descents and as many ascents. Let us describe one descent.
\begin{itemize}
\item The input is down-sampled by a convolution of kernel size 2, stride 2, and that doubles the depth of its input. The stride allows to halve the length of the input, incidentally blurring the precise locations of the highlighted features.
\item Then are applied two successive convolutions of kernel size $p$ and that conserve the depth of their input, both followed by the softplus function.
\item The output of this double convolution is added to its input. This "shortcut" is represented by the dotted arrows in Figure \ref{fig:v-net}. This way the convolutions produce additive (or residual) modifications. This is the principle of the famous Resnet. It has been shown that it accelerates the training process and limits the problems of vanishing or exploding gradient \cite{he}. 
\end{itemize}
The output of this descent is a signal with half the length and double the depth of the input signal.

\item [Ascending part] It works as a mirror of the descending part. Each ascent simultaneously increases the length of the signal and decreases its depth by a factor of 2. Technically, we up-sample the data using a transposed-convolution \cite{dumoulin} of size 2 followed by two convolutions of kernel size $p$ with softplus activations.

\item [Long shortcuts] At each level, the features obtained in the ascents are combined with the features created by the descents, using a summation. This allows to retrieve some fine resolution details, which are lost by the blurring produced by down and up samplings.  

\item [Weighted mean] The resulting signal of depth $d$ is turned into a one dimensional output $Y$ by simply taking a weighted mean of all $d$ channels. This operation is implemented with a convolution of kernel size 1. No activation function is added, as we deal with a regression problem.
\end{description}
All the coefficients of the convolutions constitute the set of parameters $\theta$ of the neural network.

Table \ref{tab:parameters} summarizes the hyper-parameters we choose in practice as a reference. This choice was mostly motivated by the results given in Section \ref{sec:res-hyper-param}.
\begin{table}[h]
\begin{center}
    \begin{tabular}{lr}
        \toprule
        Hyper-parameter & Value \\
        \midrule
        size of the input window ($N$) & 512 \\
        number of levels ($\ell$) &  5\\
        depth ($d$) & 4 \\
        size of the kernels ($p$) & 11 \\
        activation function & softplus \\
        \bottomrule
    \end{tabular}
\end{center}
\caption{Hyper-parameters of the reference neural network}
\label{tab:parameters}
\end{table}
In order to avoid the oscillations that can be produced by the upsampling operation of the V-Net architecture, we initialize the transposed convolutions with constant kernels.

\subsection{Training of the neural network}
\label{sec:training}

The training consists in finding an optimal set of parameters $\hat \theta$ for the neural network. It is defined as the minimizer of a loss function that measures the error of the network: 
\begin{equation*}
\hat \theta = \underset{\theta \in \Theta}{\text{argmin}}\ \text{Loss}(\theta).
\end{equation*}
For this loss function we choose the mean absolute error (MAE) between the output of the network and the expected heat flux over a \emph{training dataset}:
$$
\text{Loss}(\theta) = \sum_{j \in \text{training dataset}} | \text{NN}_\theta (\tilde X_j) - \tilde Y_j |,
$$
where $(\tilde X_j)$ (resp. $\tilde Y_j$) are data of size $N \times 4$ (resp. $N$), obtained from kinetic simulations after pre-processing (P) and slicing (Sl) of vectors $(\eps, \rho, u, T)$ (resp. $q$).
In practice, the minimization is carried out with a gradient descent algorithm.

This whole process falls under the category of supervised learning, as it requires a \emph{labelled} dataset, including both the inputs and the corresponding expected outputs.
The first step to train the neural network is thus to generate such a labelled dataset. To do so, we consider the kinetic model \eqref{eq:vlasov}-\eqref{eq:poisson}, that describes the particles' distribution from which can be derived their density, mean velocity, temperature, and heat flux. We run many simulations with a given time step and given number of discretization points in space $N_x$ and velocity $N_v$, with the numerical method described in appendix \ref{sec:vlasov-num}. We save the results in our dataset at selected times.

The making of this dataset is described in more details in Section \ref{sec:data}. In particular, it requires crucial choices in the initial distributions, the recording times and the range of the parameter $\eps$. It has a direct impact on the range of validity of the closure obtained and it will be  numerically discussed in Section \ref{sec:results}.

In the end, the network is trained on a dataset of $10\,000$ entries $(X_k ; Y_k) = (\varepsilon_k, \rho_k, u_k, T_k ; q_k)$, that are processed as described in Section \ref{sec:processing}. In particular, each vector (of size $N_x=1\,024$) is sliced into overlapping windows of size $N=512$, resulting in $80\,000$ inputs and labels. 4\% of these are isolated as a small test dataset and 10\% of what remains are used for validation. The rest is used for the actual training.

This training consists in 5 series of 120 epochs with a decaying learning rate, reset to its initial value of $0.005$ between each series. We use the mini-batch gradient descent method, with the Adam optimizer and mini-batches of size $1\,024$, to minimize the mean absolute error (MAE).


\section{Data generation}
\label{sec:data}

In order to train a neural network to predict the heat flux $q$ knowing the Knudsen number $\varepsilon$, the density $\rho$, the mean velocity $u$ and the temperature $T$, we need to generate a \emph{training dataset}. This dataset must meet two criteria. First, it must be a \emph{labelled} dataset, \emph{i.e.} a dataset with both the input $X=(\varepsilon, \rho, u, T)$ and the expected output $Y=q$. The expected output is the "correct" output of the given input, from which we want the network to interpolate, or \emph{generalize}, for new inputs. Then, this dataset must contain as much diversity as possible, for the interpolation to be usable in many different situations.

\begin{figure}
    \centering
    \includegraphics[scale=1]{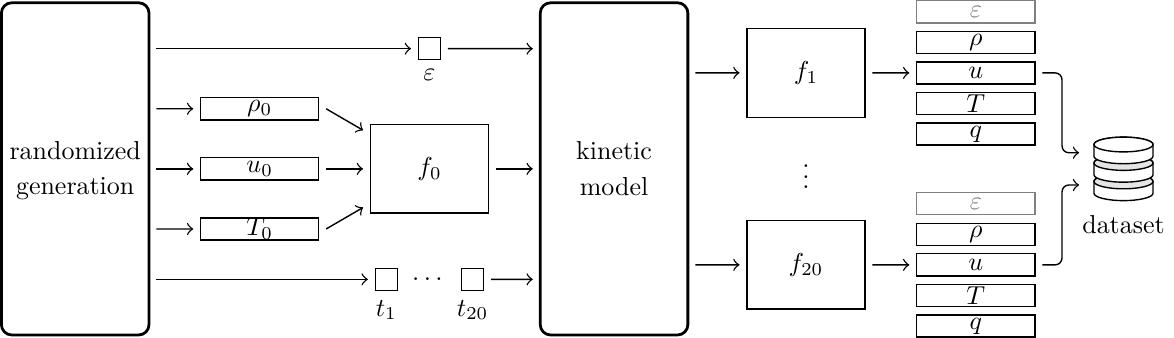}
    \caption{Scheme of the data generation process for one simulation with the kinetic model. With one Knudsen number and one initial solution, we produce 20 entries in the dataset.}
    \label{fig:generation}
\end{figure}

Consequently, the generation of the training dataset relies on two key ingredients. First, a solver for the kinetic model that can reliably estimate the distribution of particles in the phase space, from which can be derived all the fluid quantities $\rho$, $u$, $T$ and $q$. Second, a mechanism to produce a variety of distributions that can be given to the solver as initial solutions. This whole process of data generation is illustrated in Figure \ref{fig:generation} and described in more details in this section ; except for the solver for the kinetic model that is not specific to our work and is described in appendix \ref{sec:vlasov-num}.

\subsection{Global description}

Let us first get a glance at the global process illustrated by Figure \ref{fig:generation}. To begin with, a randomized process allows us to produce a Knudsen number $\varepsilon$, an initial solution $f_0$, and a \emph{recording time} $t$. These three elements are given to the kinetic model that computes the evolution of the distribution up to $t$. The resulting distribution $f$ is then used to compute the density $\rho$, the mean velocity $u$, the temperature $T$ and the heat flux $q$ of the system. On top of these four quantities, we add the Knudsen number $\varepsilon$ as a constant vector of the same size as the other four, so that it can be used as input by the network. This could also allow us to generalize our method to a problem with a Knudsen number that varies in space. In the end, these five vectors $(\varepsilon, \rho, u, T ; q)$ form an \emph{entry} in the dataset.

Actually, in order to capitalize on the simulations that can be computationally expensive —especially in higher dimensions—, we decide to produce several entries with each simulation. For a given Knudsen number and a given initial solution, we generate not one but 20 different recording times $t_1 < \cdots < t_{20}$, and compute the distribution at each one of these times, respectively $f_1$, ..., $f_{20}$. From each one of these distributions can be derived $\rho$, $u$, $T$, $q$, to which we prepend $\varepsilon$ to form 20 different entries to be stored in the dataset. Thus, for the cost of one simulation up to $t_{20}$, we produce 20 different entries corresponding to different intermediate steps in the simulation.

To build a whole dataset, we repeat this process with 100 different values of $\varepsilon$, and with 5 different randomly selected initial solutions for each one of these $\varepsilon$. Thus, we get a dataset with $10\,000$ entries, based on $500$ different initial solutions. We produce two such datasets: a training dataset used to train the neural networks, and a test dataset used to measure and compare the performance of the trained networks on new data.

\subsection{Random initial conditions}
\label{sec:solinit}

\begin{figure}
    \centering
    \includegraphics[scale=1]{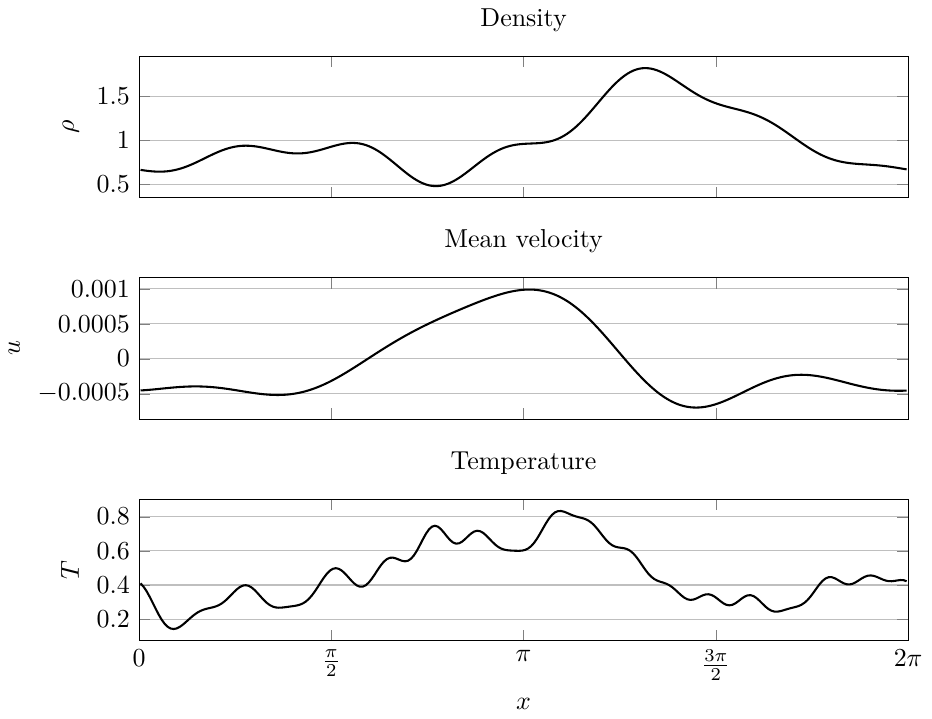}
    \caption{Possible functions $\rho$, $u$ and $T$ generated with the process described in Section \ref{sec:solinit}.}
    \label{fig:solinit}
\end{figure}

The initial solution $f_0$ is a Maxwellian as in (\ref{eq:maxwellian}), with density $\rho$, mean velocity $u$ and temperature $T$ randomly generated as partial Fourier series under the form
$$
\alpha\times\left(\frac{a_0}{2} + 0.5\sum_{n=1}^N (a_n \cos(nx)+b_n\sin(nx))\right), \quad x \in [0, 2\pi].
$$
We set $N$ to 20 and we randomly generate the $a_n$ and $b_n$ coefficients for $n\geqslant1$ according to a uniform distribution on $\left[-\frac{1}{n}, \frac{1}{n}\right]$. This choice for the $a_n$ and $b_n$ coefficients set the high frequencies to low amplitudes, resulting in more regular functions. Then, the choice for $\alpha$ and $a_0$ depends on the function to be generated. For the density, $\frac{a_0}{2}=1$ and $\alpha=1$. For the temperature, $\frac{a_0}{2}=1$ and $\alpha$ is chosen uniformly in $[0.1, 1]$. For the mean velocity, $\frac{a_0}{2}$ is chosen uniformly in $[-1, 1]$, and $\alpha$ is chosen so that the maximum Mach number
$$\underset{x \in [0,2\pi]}{\max} \frac{|u(x)|}{\sqrt{2T(x)}},$$
falls uniformly in $[10^{-4}, 5\cdot 10^{-1}]$ in logarithmic scale. We also make sure that the density and temperature generated are positive functions.
Possible functions generated with this process are given Figure \ref{fig:solinit}.
For their discrete form, these functions are sampled with $N_x = 1\,024$ points to get vectors, that are turned into a discrete Maxwellian our solver can work with, with $N_v=141$ points is velocity.

\subsection{Random Knudsen numbers and recording times}

In this work we focus on Knudsen numbers in the range $[\varepsilon_{\min}, \varepsilon_{\max}] = [0.01, 1]$. First because below $0.01$ the fluid model with the Navier-Stokes estimation is already very accurate, and also because we do not want a range too wide.
Regarding the distribution of the values in this interval, we want a decent amount of entries with $\varepsilon$ really close to $0.01$, but we do not want too few entries for $\varepsilon$ close to $1$. As a result, a uniform or logarithmic distribution over the interval $[0.01, 1]$ is not satisfying. Instead we opt for a uniform distribution for $\sqrt{\varepsilon}$ in the interval $[\sqrt{\varepsilon_{\min}}, \sqrt{\varepsilon_{\max}}]$. This allows us to have about one fourth of the Knudsen numbers between $0.01$ and $0.1$, and three fourths between $0.1$ and $1$.
Also we use a deterministic distribution for the train dataset in order to have a nice coverage of the interval, while we use a random distribution for the test dataset in order to test the ability of the trained networks to generalize to new Knudsen numbers.

For the recording times $t_1 < \cdots < t_{20}$, we do not want $t_1$ to be too close to zero as the initial maxwellian state always has a heat flux of zero and does not carry much physical information, and we do not want $t_{20}$ to be too big so that the simulations do not last too long. For these reasons we choose $20$ times uniformly in the interval $[0.1, 2]$, that are then sorted and given successively to the solver for the kinetic model.

\subsection{Computing of the fluid quantities}

Once the Knudsen number, the initial solution and the times are generated, a simulation is computed with the numerical method described in appendix \ref{sec:vlasov-num}. It relies on a discretization of the phase space
\[
(x_i, v_j)\quad,\quad i\in\{1,...,N_x\},\  j\in\{1,...,N_v\},
\]
and results in an approximation $(f_{i,j})$ of the real distribution $f$:
\[
f_{i,j} \simeq f(x_i, v_j).
\]
From this distribution can be derived the density $\rho$, mean velocity $u$, temperature $T$ and heat flux $q$ needed for the dataset, as mentioned in Section \ref{sec:kinetic}. In their discrete form, these quantities are vectors computed as follows: for any $1 \leqslant i \leqslant N_x$,
\[
\rho_i = \sum_{j=1}^{N_v}f_{i,j}, \quad \rho_i u_i = \sum_{j=1}^{N_v} v_jf_{i,j}, \quad \rho_i T_i = \sum_{j=1}^{N_v} (v_j-u_i)^2 f_{i,j}
\]
and
\[
\quad q_i = \sum_{j=1}^{N_v} (v_j-u_i)^3 f_{i,j}.
\]
For our datasets, we use a discretization with $N_x=1\,024$ points in space, and $N_v=141$ points in velocity.

\section{Data processing}
\label{sec:processing}

\begin{figure}
    \centering
    \includegraphics[scale=1]{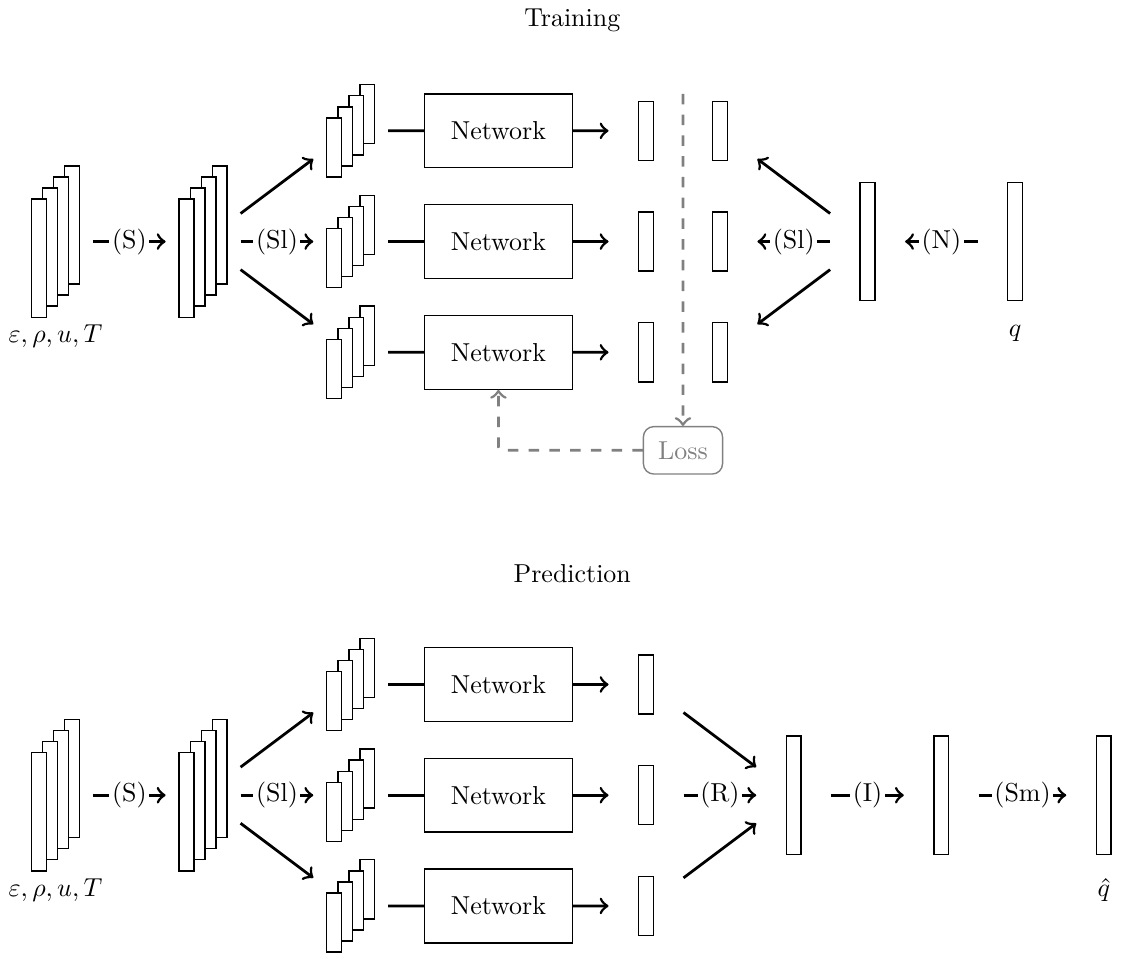}
    \caption{Scheme of the whole data processing for the training and the prediction processes respectivily. The different operations are (S) standardization, (Sl) slicing, (N) normalization, (I) inverse normalization, (R) reconstruction and (Sm) smoothing.}
    \label{fig:processing}
\end{figure}

In this section we describe how the data is processed on both sides of the neural network, both for the training process and the prediction process. At this point we have generated some data with the kinetic model, and want to process it for a neural network that takes as input a 1D signal with four \emph{channels} —the Knudsen number, the density, the mean velocity and the temperature—, and outputs a 1D signal with one channel —the heat flux. This data processing consists in the transformations described below and summarized in Figure \ref{fig:processing}. In the remainder of this section, we use the word \emph{standardized} to describe a quantity with mean zero and unit variance, and \emph{normalized} a quantity with norm one or that has been applied the transformation described in Section \ref{sec:out-norm}.

\subsection{Input standardization}
\label{sec:in-std}

Neural networks are known to be easier to train with standardized inputs, where all channels have the same order of magnitude. Consequently, we standardize the inputs channel by channel, by removing the mean and scaling to unit variance the whole training dataset. For example for the density $\rho$, if $( \rho^k )$ is the family of densities of the whole training dataset, this standardization applied on a given entry $\rho^{k_0}$ can be written
\[
\rho^{k_0}_\text{standard} =
\frac{\rho^{k_0} - \underset{k,i}{\mathrm{mean}}\ \rho^k_i}
{\underset{k,i}{\mathrm{std}}\ \rho^k_i},
\]
where \emph{mean} is the empirical mean and \emph{std} the empirical standard deviation.

Since the neural network is trained with this standardized data, it is designed to work only with inputs that received that specific standardization: after the training, any new input given to the neural network must be applied that very same standardization. In particular, the means and standard deviation computed in the training dataset have to be stored with the neural network, in order to be available later to make predictions for the fluid model.

\subsection{Output normalization}
\label{sec:out-norm}

 To prevent the outputs in our datasets from being smaller than the typical error of our neural network, we choose to train the network to predict normalized output. Otherwise, the network would produce pure noise when trying to predict them, which turns to be problematic when used by the fluid model, especially regarding its stability. Another solution would be to use a cost fonction measuring a relative error instead of an absolute one but it turns out to be less efficient. Let us describe how we proceed.

The main issue when training the network with normalized output is the reversibility of the normalization. Since in the end we want to predict a heat flux, and not a normalized heat flux, we need to be able to apply an inverse normalization at the output of the neural network. As a consequence, this transformation cannot use any knowledge on the heat flux to be predicted, that would be available in the labelled training dataset, but would not in a real case scenario. For this reason, we choose to normalize the heat flux with its estimation given by the Navier-Stokes approximation (\ref{eq:ns}), and that can be computed from the Knudsen number, the density and the temperature, all available in a real case scenario.

But using an estimation of the heat flux brings another concern: when this estimation is way off, the normalization (or inverse normalization) would distort the information a lot. And this is expected to happen, as the Navier-Stokes approximation tends to greatly overestimate big heat fluxes with Knudsen numbers in the range we consider. To prevent this from happening, we choose to normalize heat fluxes only when the corresponding Navier-Stokes estimation is below a given threshold $\theta$, that we set to $0.1$. This threshold is also a way of only addressing the heat fluxes that were problematic in the first place.

Thus, for a given input $(\varepsilon^{k_0}, \rho^{k_0}, u^{k_0}, T^{k_0})$ in the training dataset with expected output $q^{k_0}$, the normalization we use can be written
\[
q^{k_0}_{\text{norm}} = \left\{
\begin{array}{cl}
    \displaystyle{q^{k_0}\times\frac{\theta}{q^{k_0}_{NS}}}, & \text{if } 0 < q^{k_0}_{NS} \leqslant \theta, \\
    &\\
    q^{k_0}, & \text{otherwise}, \\
\end{array}\right.
\]
with
\[
q^{k_0}_{NS} = \max_{i=1,\ldots , N} \left| \frac{3}{2}\varepsilon^{k_0}\rho^{k_0}_i(\partial_x T)^{k_0}_i \right|,
\]
where $i$ refers to the index of the vectors and the spatial derivative is approximated with a centered finite difference formula.
Note that we choose not to normalize at all heat fluxes such that $q_{NS} = 0$.
Once trained to predict these normalized heat fluxes, the neural network must be followed by an inverse normalization when used to make predictions of non normalized heat fluxes. This inverse normalization simply reads
\[
q^{k_0} = \left\{
\begin{array}{cl}
    \displaystyle{q^{k_0}_{\text{norm}}\times\frac{q^{k_0}_{NS}}{\theta}}, & \text{if } 0 < q^{k_0}_{NS} \leqslant \theta, \\
    &\\
    q^{k_0}_{\text{norm}}, & \text{otherwise}. \\
\end{array}\right..
\]

\subsection{Slicing and reconstructing}
\label{sec:slicing}

In order to be usable with meshes of different sizes, the neural network is designed to work with only one portion of the signal at a time as input, and to return the corresponding portion as output. As a consequence, each signal given to the network has first to be sliced into several \emph{windows} of the right size. For the training process, both the input and the expected output have to be applied the same slicing. For the prediction process, the input has to be sliced into windows, and the resulting predictions have then to be aggregated to reconstruct a complete signal.

In one case or the other, we slice the signal into overlapping windows. The overlapping serves two purposes. For the training process, it allows us to artificially increase the amount of data we can give to the network from the training dataset. Though it introduces redundancy, such data augmentation is known to help prevent overfitting and improve the accuracy of the network \cite{data-augmentation}. For the prediction process, we use this redundancy introduced by the overlapping when reconstructing the output signal, to provide a first smoothing process and to dismiss the defects on the borders of the predictions.

\begin{figure}[t]
    \centering
    \includegraphics[scale=1]{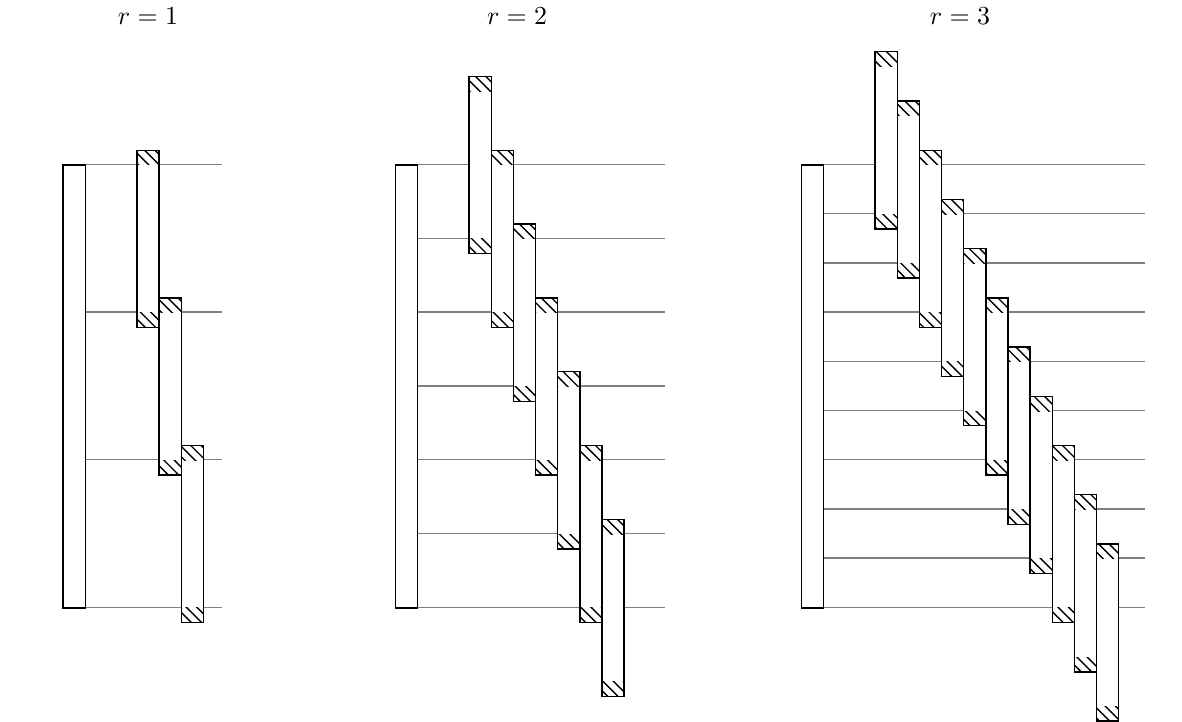}
    \caption{Slicing of a 1D signal into overlapping windows with different redundancy parameters $r$. The original signal is represented vertically on the left, and the windows are on the right. The hatched parts of each window are the margins that are ignored when reconstructing the signal. The windows that exceed the span of the original signal are completed by periodicity.}
    \label{fig:slicing}
\end{figure}

These defects come from the necessary padding in our V-Net architecture, that can cause oscillations of big amplitudes at the ends of the predictions. As a consequence, each window must have a \emph{margin} that will be ignored when reconstructing the signal. We set this margin to 10\% of the output on each side, leaving 80\% actually useful for the reconstruction, later referred to as the \emph{useful part}. This observation on its own is enough to require some overlapping between the windows to reconstruct the whole signal, but we introduce even more overlapping. We slice the original signal such that each one of its points is found in a fixed number $r$ of windows, margins aside. We call this number the \emph{redundancy parameter}. The bigger it is, the bigger the overlapping, and the more windows are produced. This slicing is illustrated Figure \ref{fig:slicing}.

\begin{figure}
    \centering
    \includegraphics[scale=1]{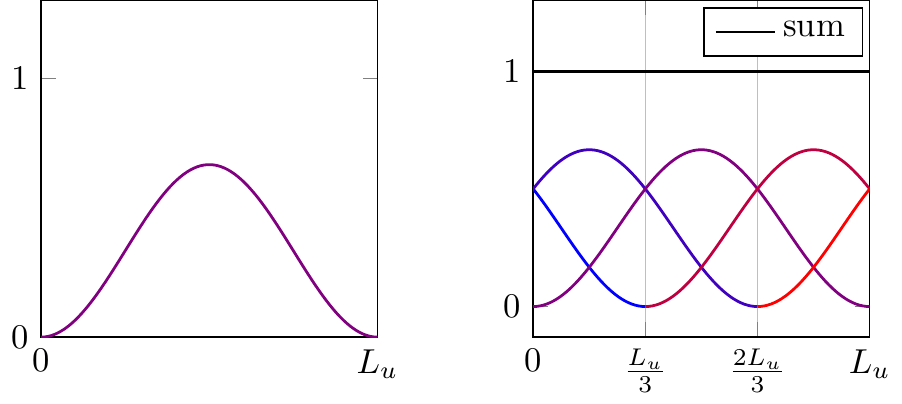}
    \caption{Kernel used for the reconstruction of an entire signal from the windows predicted by the neural network, with the redundancy parameter $r=3$.}
    \label{fig:reconstruction}
\end{figure}

For the reconstruction of an entire signal from the predictions of the network, we get rid of the margins and multiply the useful parts by the kernel
\[
\frac{1}{n}\left(\cos\left(\frac{2\pi}{L_u}x-\pi\right) + 1\right), \quad x\in[0,L_u]
\]
where $r$ is the redundancy parameter and $L_u$ the length of the useful part in the physical space. We then sum up the resulting windows.
This kernel has the property to sum up to one when duplicated every $\frac{L_u}{r}$ (see Figure \ref{fig:reconstruction}). As a consequence, each value in the reconstructed signal is a weighted mean of the $r$ values predicted by the neural network, which results in a somewhat smooth estimation of the heat flux.

\subsection{Smoothing of the output}

At this point we have a process that allows us to estimate the heat flux $q$ with the neural network, and that can be used in a fluid model. But as it is, even with the smoothing provided by the aggregation of different predictions as described in the previous section, this estimation can show some small irregularities. In the fluid model, these irregularities on $q$ are amplified by the computation of $\partial_x q$, and can cause the instability of the hyperbolic system. This instability cannot be prevented by refining the time discretization, we have to control the irregularities in the estimated heat flux. To do so, we smooth the reconstruted heat flux thanks to a convolution with a gaussian kernel. This kernel with standard deviation $\sigma$ reads
\[
w:t\in[-3\sigma, 3\sigma] \mapsto \frac{e^{-\frac{1}{2}\frac{t^2}{\sigma^2}}}{\int_{-3\sigma}^{3\sigma} e^{-\frac{1}{2}\frac{t'^2}{\sigma^2}}\,dt'},
\]
and the smoothed heat flux is
\[
\Tilde{q}(x) = \int_{-3\sigma}^{3\sigma} q(x+t)w(t)\,dt.
\]
We show in Section \ref{sec:res-smoothing} that $\sigma$ can be chosen relatively large ($\sigma \simeq 0.05$) without impacting the accuracy of the estimation, and has to be for the method to be stable.


\section{Results}
\label{sec:results}

This section is divided in three parts: first we introduce the results regarding the neural network only, then when it is used in the fluid model, and finally how the whole model performs in configurations that differ from the one used for the training. Unless otherwise stated, the hyperparameters of the neural network and the training are those presented in section \ref{sec:architecture}.

\subsection{Accuracy of the network}
\label{sec:accuracy}

In this section we show numerical results of the neural network, independently of its use with the fluid model. To this end we use the test dataset with entries $(\varepsilon,\rho,u,T;q)$ computed with the kinetic model, or directly compare the predictions with a simulation of the kinetic model. To measure the accuracy of the network on a given entry, we use the relative error in norm $L^2$
\[
\frac{|| q - \hat{q} ||_2}{|| q ||_2},
\]
where $\hat{q}$ is the prediction of the neural network with inputs $(\varepsilon,\rho,u,T)$. We also observed the relative error in norm $L^\infty$, but it was essentially the same, that is why we focus on the $L^2$ norm in what follows.

Unless stated otherwise, the results given use a redundancy parameter of 2 and a smoothing of $\sigma \simeq 0.06$. This last choice is motivated by the results of Section \ref{sec:res-smoothing}.

\subsubsection{Comparison with the Navier-Stokes estimation}

\begin{figure}
    \centering
    \includegraphics[scale=1]{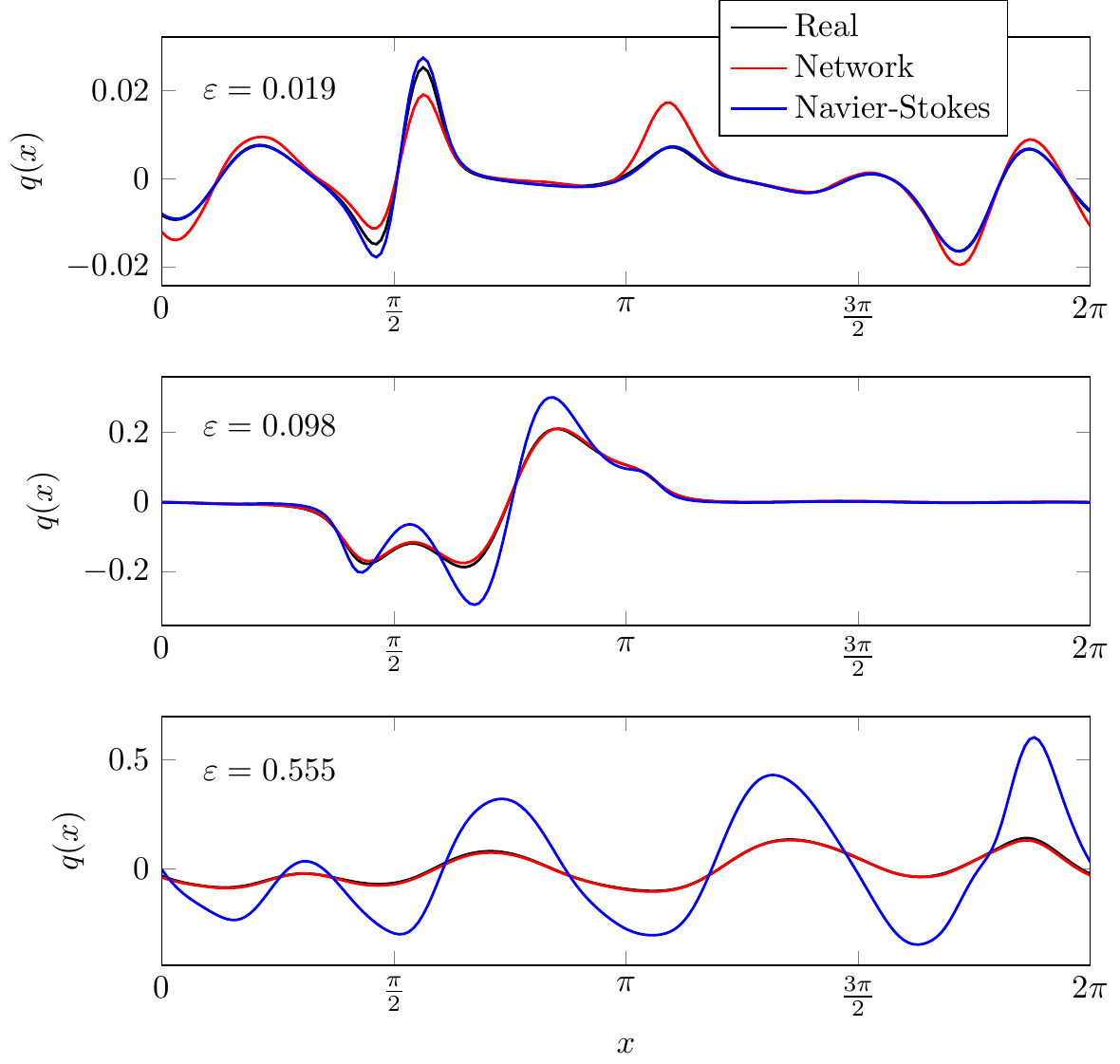}
    \caption{Three examples of predictions.}
    \label{fig:pred-ex}
\end{figure}

Figure \ref{fig:pred-ex} shows some examples of predictions, compared to the real heat flux and the estimation of Navier-Stokes (\ref{eq:ns}). The first example illustrates how the Navier-Stokes estimation tends to perform better on data with small Knudsen numbers. The second and third example show how on the contrary it very often overestimate the heat flux when the Knudsen number increases.

\begin{figure}
    \centering
    \includegraphics[scale=1]{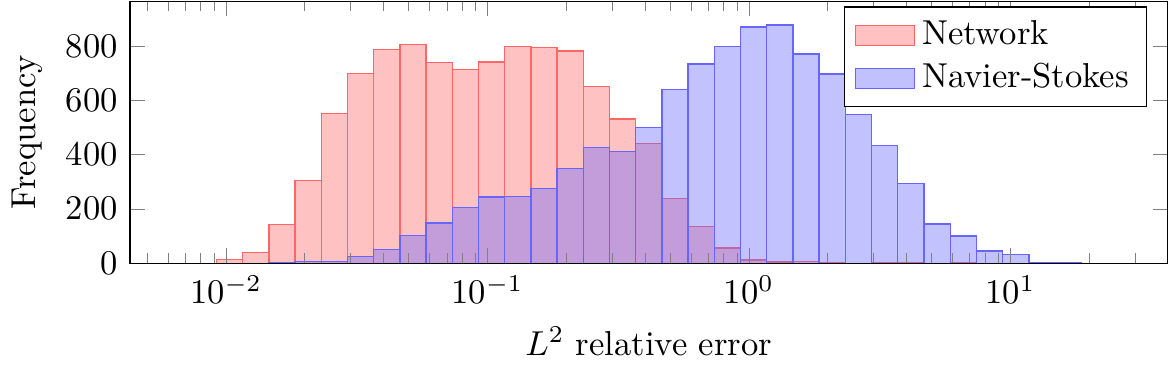}
    \caption{Distribution of the relative errors of the neural network and the Navier-Stokes estimation over the test dataset ($10,000$ predictions).}
    \label{fig:pred-err}
\end{figure}

The relative errors of both the network and the Navier-Stokes estimation over the whole test dataset are summarized in the histogram Figure \ref{fig:pred-err}. On this dataset, the overall error of the neural network is about an order of magnitude below that of the Navier-Stokes estimation.

\begin{figure}
    \centering
    \includegraphics[scale=1]{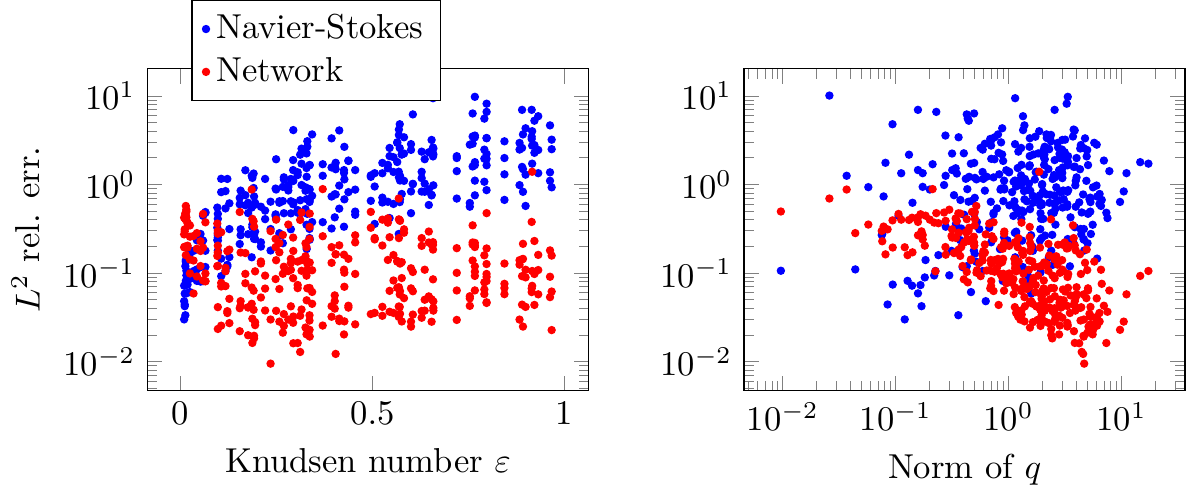}
    \caption{Error of the Navier-Stokes estimations and of the network on the test dataset, depending on the Knudsen number and the $L^2$ norm of the real heat flux. Each dot corresponds to one entry in the test dataset. For better clarity, only one $31^{\text{th}}$ of the data is shown here.}
    \label{fig:pred-eps-q}
\end{figure}

\begin{figure}
    \centering
    \includegraphics[scale=1]{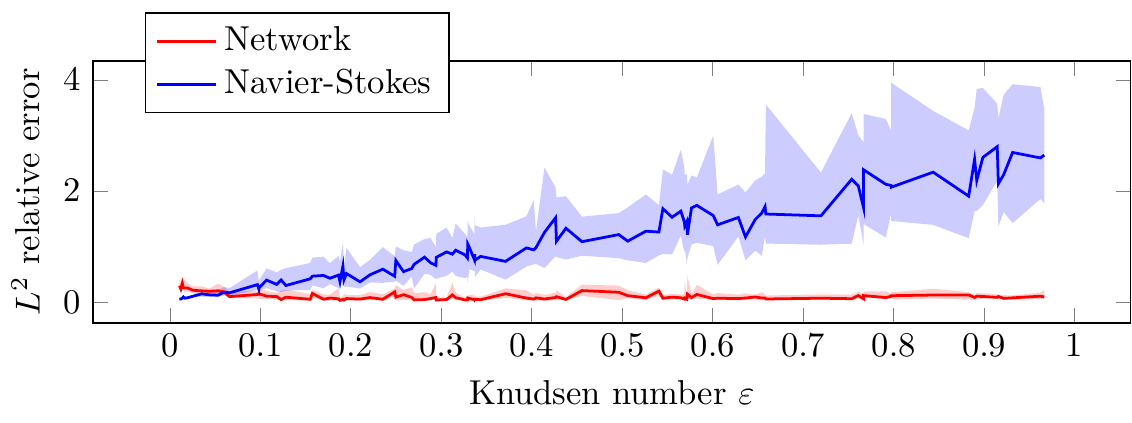}
    \caption{Relative error of the Navier-Stokes estimations and of the network on the test dataset, depending on the Knudsen number. The line represents the median of the error over the 100 entries of the test dataset for each epsilon, and the coloured area the interquartile interval.}
    \label{fig:pred-eps}
\end{figure}

Figure \ref{fig:pred-eps-q} gives a closer look at these errors and highlights two important factors in the performance of the network over the Navier-Stokes estimation: the Knudsen number and the norm of the real flux. The first scatter plot (Fig. \ref{fig:pred-eps-q}, left) shows that the error of the Navier-Stokes estimation heavily depends on the Knudsen number $\varepsilon$, which is not the case for the neural network predictions, except for small values. This is even better illustrated by Figure \ref{fig:pred-eps}, that uses the same data but where the error is not set in logarithmic scale.
The second scatter plot (Fig. \ref{fig:pred-eps-q}, right) shows that there is a strong correlation between the norm of the real heat flux and the relative error of the network: the smaller the heat flux, the bigger the error of the network. This could be explained by a lack of small heat fluxes in the training dataset or by the normalization we chose for small heat fluxes. In any case, it should not be that big of a problem when used with the fluid model, as smaller heat fluxes have a smaller impact on the result, at least as long as they are quite smooth.

\begin{figure}
    \centering
    \includegraphics[scale=1]{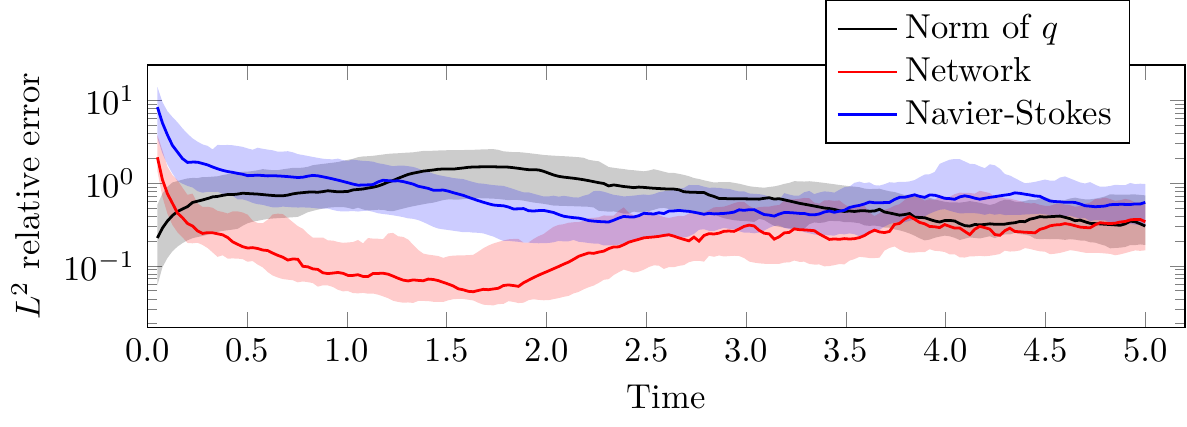}
    \caption{Norm of the real heat flux and relative errors of the predicted heat flux and the Navier-Stokes estimation throughout simulations. Median and interquartile interval over 50 simulations.}
    \label{fig:pred-simu-err}
\end{figure}

It can also be interesting to look at the evolution of the error during a simulation, to get a better idea of what will happen when used with the fluid model. To do so, we choose 50 random initial solutions and 50 random Knudsen numbers in the same way we did to build the test dataset, and run the kinetic model with these data. At regular times, the heat flux of the kinetic model is compared to the heat flux the network would have predicted and to the Navier-Stokes estimation. Figure \ref{fig:pred-simu-err} shows the statistical results of this experiment. We observe that the network has a significant advantage over the Navier-Stokes estimation during the first $2.5$ time units, before decreasing to a smaller advantage. This result can be explained by the fact that the heat flux is bigger during the first three time units (the dissipation seems decrease the size of $q$ in time), and by the strong correlation between the norm of the heat flux and the error of the network, also visible in this figure.

\subsubsection{Hyper-parameters of the neural network}
\label{sec:res-hyper-param}

\begin{figure}
    \centering
    \includegraphics[scale=1]{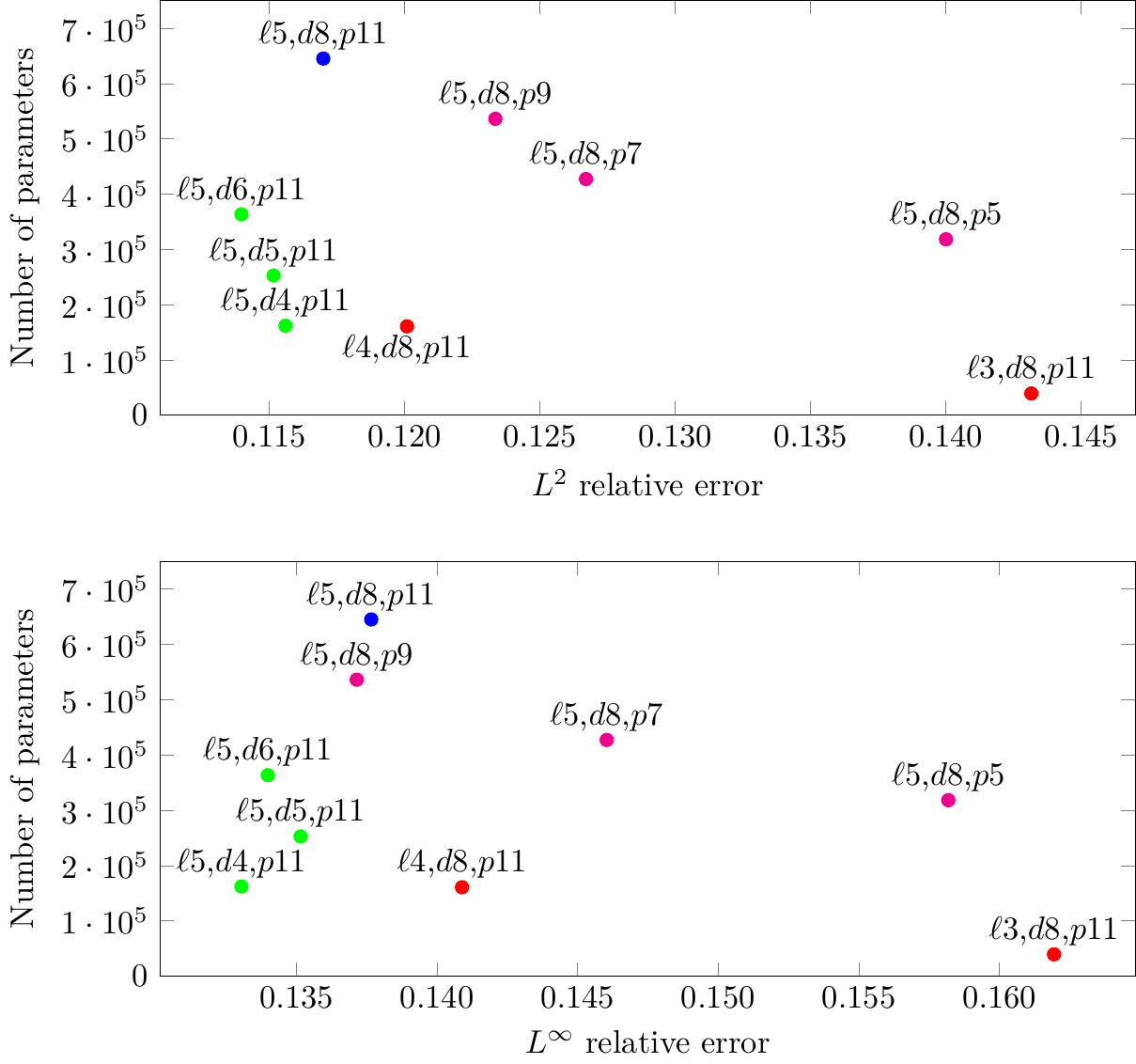}
    \caption{Median relative error and number of learnable parameters for V-Nets with different hyper-parameters : ($\ell$) number of levels, ($d$) depth, and ($p$) size of the kernels. Points sharing the same color represent networks that differ by only one hyper-parameter.}
    \label{fig:networks-err}
\end{figure}

In this section we show the results that motivated the choice for some hyper-parameters of the V-Net. We focus on the number $\ell$ of levels, the depth $d$, and the size $p$ of the kernels of the convolutions. Figure \ref{fig:networks-err} compares the median relative error of V-Nets with different sets of hyper-parameters. For these results the window size was set to 256, but we later decided to set it to 512 as it gave slightly better results.

We observe that all hyper-parameters do not affect the performance of the V-Net in the same way. For instance, decreasing the size of the kernels significantly decreases the accuracy of the trained network, while not decreasing the number of parameters as much. On the other hand, decreasing the depth allows to decrease the number of parameters with a very small effect on the accuracy.
In the remainder of this paper we work with the "l5,d4,p11" architecture, as it seems to be a good compromise between the number of parameters and the accuracy of the network.

\subsection{Fluid model with the neural network}
\label{sec:res-fluid}

In this section we look at how the neural network performs when used in the fluid model. We denote by "Fluid+Network" this method. The fluid model is solved using an explicit finite volume method as presented in appendix \ref{sec:fluid-num-explicit}. We compare it to three others:
\begin{description}
    \item [Kinetic]: the kinetic model. It is the most accurate model and serves as a reference for the real target. The numerical method is described in appendix \ref{sec:vlasov-num}.
    \item [Fluid+Kinetic]: the fluid model with the heat flux from the kinetic model. It is the result we would get with a perfect neural network that makes no error, and helps to distinguish the error of the fluid model from the error of the neural network on the heat flux. The numerical method is identical to the one of the Fluid+Network method (see appendix \ref{sec:fluid-num-explicit}). Theoretically this model should give the same result as the kinetic model, but since the numerical schemes and viscosities are different it is not always true in practice.
    \item [Navier-Stokes]: the fluid model with the Navier-Stokes estimation. It does not use the same numerical method as our model, since the formula for the heat flux requires an implicit scheme to avoid too stringent stability condition (see appendix \ref{sec:ns-num}).
\end{description}
Our criteria to compare the fluid models is the $L^2$ relative error on the logarithm of the electric energy $\mathcal{E}$
$$
\mathcal{E}(t) = \int_{[0,2\pi]} E(x, t)^2\,dx,
$$
compared to the kinetic model.
For all the following tests we use a discretization of $N_x=512$ points on $[0, 2\pi]$ in physical space, and $N_v=101$ points on $[-7,7]$ in velocity space (for the kinetic model). In particular, the data is resampled on $1\,024$ points at each iteration for the neural network as explained in Section \ref{sec:integration}. The simulations are computed up to $t=8$.

\subsubsection{Examples}

\begin{figure}
    \centering
    \includegraphics[scale=1]{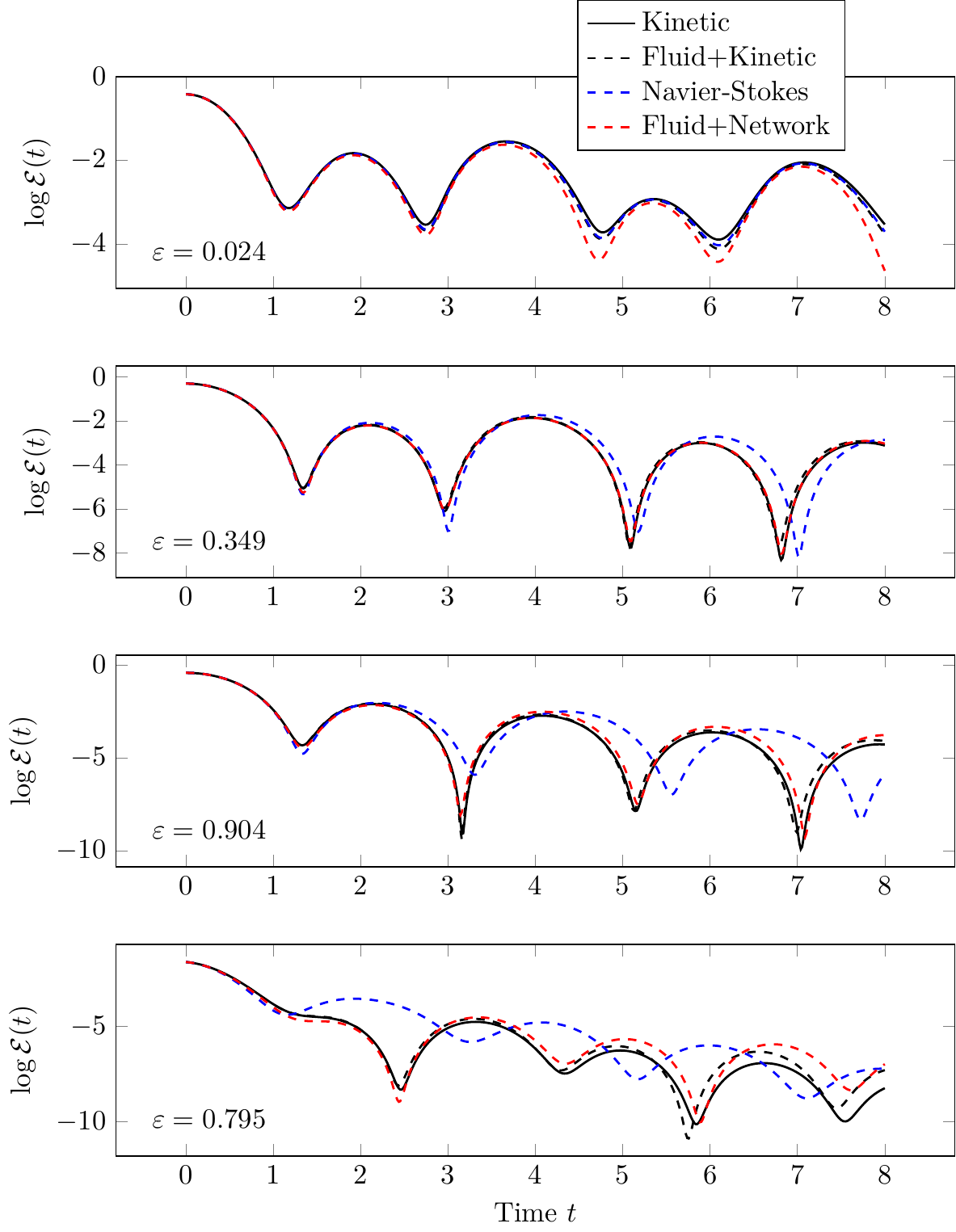}
    \caption{Examples of the evolution of the electric energy with different initial solutions and different Knudsen numbers.}
    \label{fig:res-simu-ex}
\end{figure}

Figure \ref{fig:res-simu-ex} shows examples of electric energies obtained with the four models described above. The first one uses a very small Knudsen number, and with no surprise the Navier-Stokes model performs better than the fluid model with the network. On the second example on the other hand, the Knudsen number is bigger and we observe that the oscillations are slightly shifted: Navier-Stokes shows some dispersion that the fluid model with the network does not. This dispersion often increases as the Knudsen number gets bigger as it shows on the third example. Finally the fourth one illustrates a case where the Navier-Stokes model really struggle while the fluid model with the network does pretty good.

\subsubsection{Global performances}

\begin{figure}
    \centering
    \includegraphics[scale=1]{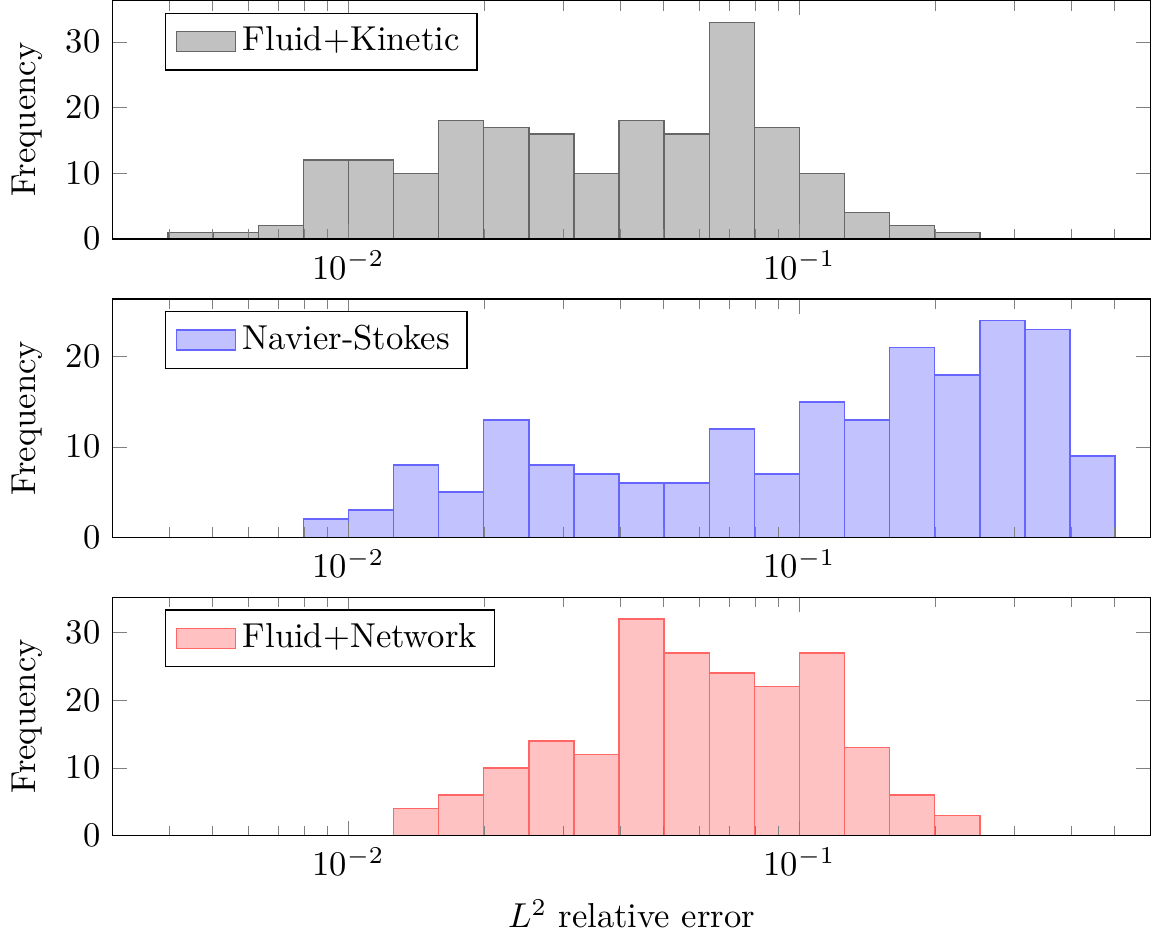}
    \caption{Distribution of the relative errors of the three fluid models on the kinetic model over 200 simulations up to $t=8$.}
    \label{fig:simu-err}
\end{figure}

The $L^2$ relative errors of the three fluid models over 200 simulations with different initial solutions and different Knudsen numbers are summarized in Figure \ref{fig:simu-err}. We observe that our closure with the neural network does not reach relative errors of $0.01$ or below, contrary to the closure with the real heat flux or even with the Navier-Stokes estimation. But it successfully maintains the relative error below $0.2$, which is close to what the closure with the real heat flux achieves, and much better than the Navier-Stokes estimation that often exceeds $0.2$. As already noted, the Fluid+Kinetic relative error should be zero but is not, due to numerical diffusion errors.

\subsubsection{Influence of the Knudsen number}

\begin{figure}
    \centering
    \includegraphics[scale=1]{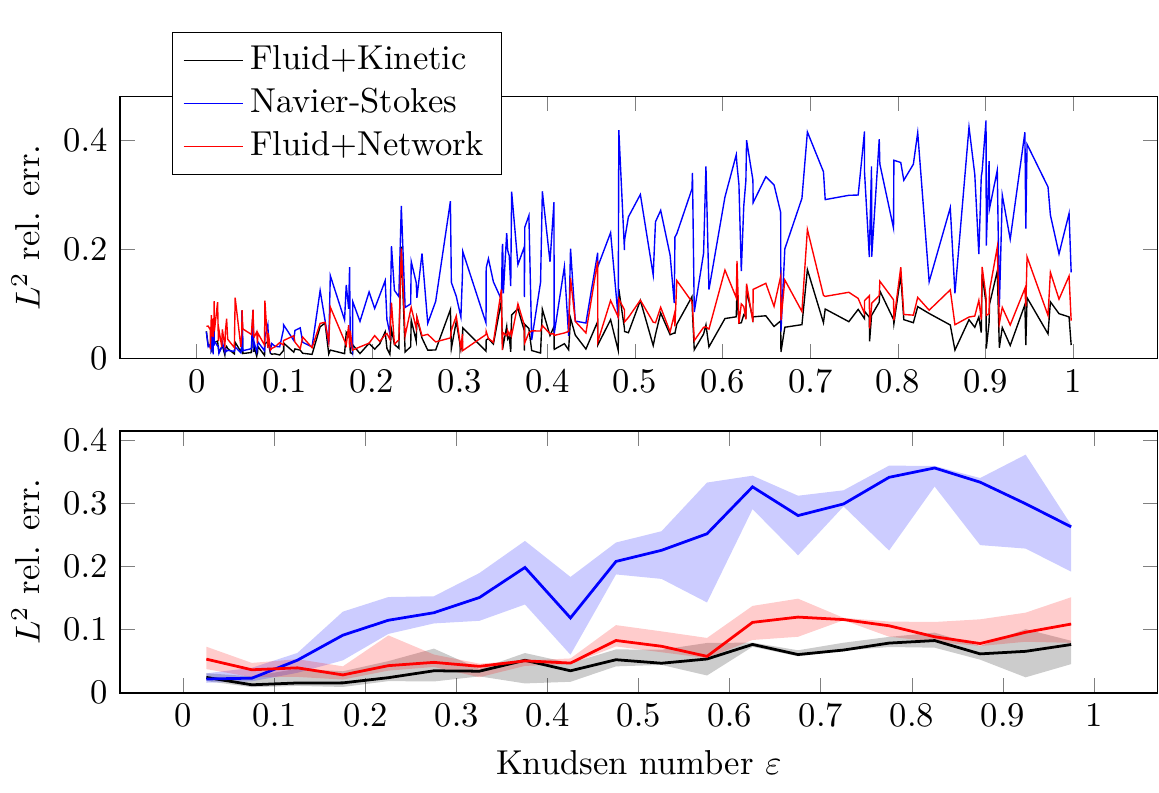}
    \caption{$L^2$ relative errors of the fluid models over the kinetic model depending on the Knudsen number. The first plot shows the raw data and the second shows the median and interquartile interval for 20 uniform classes of Knudsen numbers between $0.01$ and $1$.}
    \label{fig:res-simu-eps}
\end{figure}

In Figure \ref{fig:res-simu-eps} we look at the influence of the Knudsen number on the result of the fluid models. All three seem to decrease in accuracy as the Knudsen number increases, but it is much more significant for the Navier-Stokes model. For Knudsen numbers above $0.1$, not only the fluid model with the network seems to work systematically better than the Navier-Stokes model, but it is also quite close to the fluid model with the kinetic heat flux. This would seem to indicate that the network appropriately plays its role and that its error on the heat flux does not impact the whole model too much.

\subsubsection{Smoothing of the prediction and stability}
\label{sec:res-smoothing}

\begin{figure}
    \centering
    \includegraphics[scale=1]{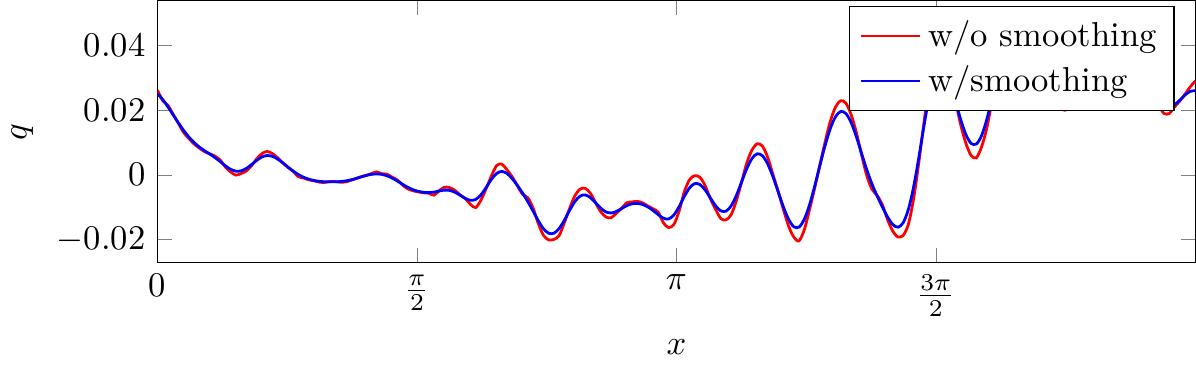}
    \caption{Example of a heat flux predicted by the neural network, with and without smoothing. The smoothing here uses $\sigma\simeq0.05$.}
    \label{fig:smooth-ex}
\end{figure}

In this section we show the impact of the smoothing on the accuracy of the predictions, as well as its impact on the stability of the fluid model using these predictions. Figure \ref{fig:smooth-ex} gives a visual example of how the smoothing modifies a prediction of the neural network.

\begin{figure}
    \centering
    \includegraphics[scale=1]{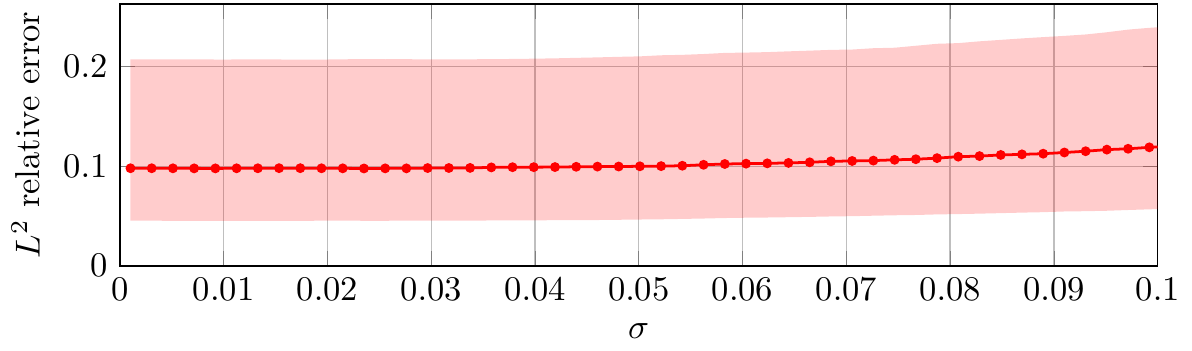}
    \caption{Relative errors of the predictions over the test dataset depending on the quantity $\sigma$ of smoothing.}
    \label{fig:smooth-accu}
\end{figure}

To plot Figure \ref{fig:smooth-accu}, the network was used to predict the heat flux of each entry in the test dataset, and different quantities of smoothing were applied before computing the relative error. We observe that the accuracy of the resulting heat flux does not deteriorate for standard deviations lower than $\sigma \simeq 0.04$, and then slightly decreases as $\sigma$ increases.

\begin{figure}
    \centering
    \includegraphics[scale=1]{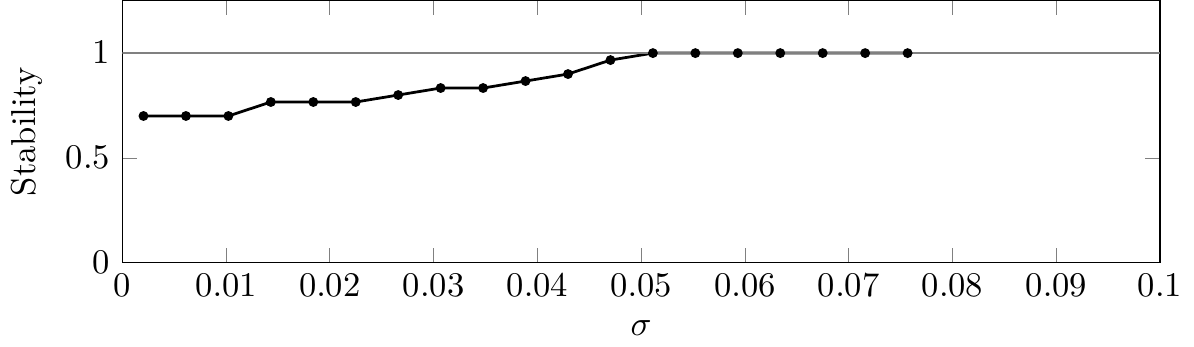}
    \caption{Proportion of simulations of the fluid model reaching $t=3$ out of 30 simulations, depending on the quantity $\sigma$ of smoothing.}
    \label{fig:smooth-stab}
\end{figure}

This smoothing was introduced to cope with the instability of the fluid model relying on the predictions of the network. To measure how this stability is improved by smoothing the output of the network, we choose 30 random initial solutions as we did to build the test dataset and 30 Knudsen numbers uniformly distributed between $0.01$ and $1$. Then we choose a set of smoothing parameters $\sigma$, and for each one of these we run the 30 simulations with the fluid model, using the network with this quantity of smoothing. Finally, for each $\sigma$ we look at the proportion of simulations that reached $t=3$, from which we assume the model will remain stable. The results are shown Figure \ref{fig:smooth-stab}. It appears that all 30 simulations reach $t=3$ for $\sigma$ above approximately $0.05$. We use $\sigma \simeq 0.06$ for our tests with the fluid model shown in the next section, and it appears to be enough since no simulation failed so far with this quantity of smoothing.

\subsection{Using the network with different configurations}

In this section we are interested in the flexibility of our approach to different configurations. First we take a look at different discretizations of space to measure the ability of the network to generalize to other resolutions than the one it has been trained with. Then we introduce discontinuities in the initial solutions to see how the fluid model with the network reacts in terms of stability and accuracy.

\subsubsection{Different resolutions}
\label{sec:res-resolution}

The network was trained with data using $1\,024$ points on $[0, 2\pi]$ in space. The fact that it uses windows of size $512$ would allow it to be used with any discretizations with at least $512$ points, or even less if we use the periodicity to make windows of size $512$. However, with initial solutions of the same form than those used to build the datasets, the change in resolution modifies the way the data is interpreted by the neural network. For instance if we sample an initial condition on $512$ points instead of $1\,024$, the signal frequencies would be multiplied by two in the eyes of the neural network. This can result in inputs differing from the training inputs, and a decrease in accuracy is to be expected.

\begin{figure}
    \centering
    \includegraphics[scale=1]{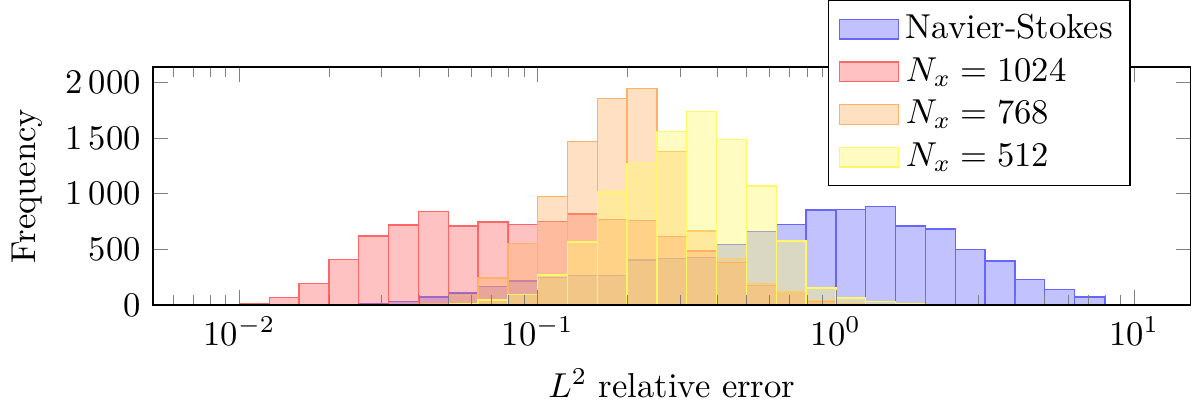}
    \caption{Distribution of the relative errors of the network on the test dataset with its original resolution, downsampled to $768$ points and downsampled to $512$ points.}
    \label{fig:pred-reso}
\end{figure}

One way to measure this effect is to evaluate the network on the test dataset but resampled at different resolutions, to mimic the data that could be generated by the fluid model with these resolutions. We can then compute the relative errors as in Section \ref{sec:accuracy}. The results Figure \ref{fig:pred-reso} show that the accuracy of the network does indeed decrease as the resolution of the data moves away from its original resolution.

To avoid this issue, if the fluid model uses a different resolution it is best to resample the data for the network to match the resolution of the training data, as this resampling operation only introduces a negligible error. However, this phenomenon shows some of the limits of the network in terms of generalization. With initial solutions that differ from the ones used to train the network —in frequency for instance—, we could expect to see a significant loss of accuracy.

\subsubsection{Discontinuities}

The datasets were built from simulations starting with a continuous initial solution, but using the network in the fluid model with discontinuous initial solutions could be an option. In practice, the initial discontinuities quickly fade away in this physical system, but it might be enough for the predictions of the network to cause some oscillations and the instability that goes with it.

\begin{figure}
    \centering
    \includegraphics[scale=1]{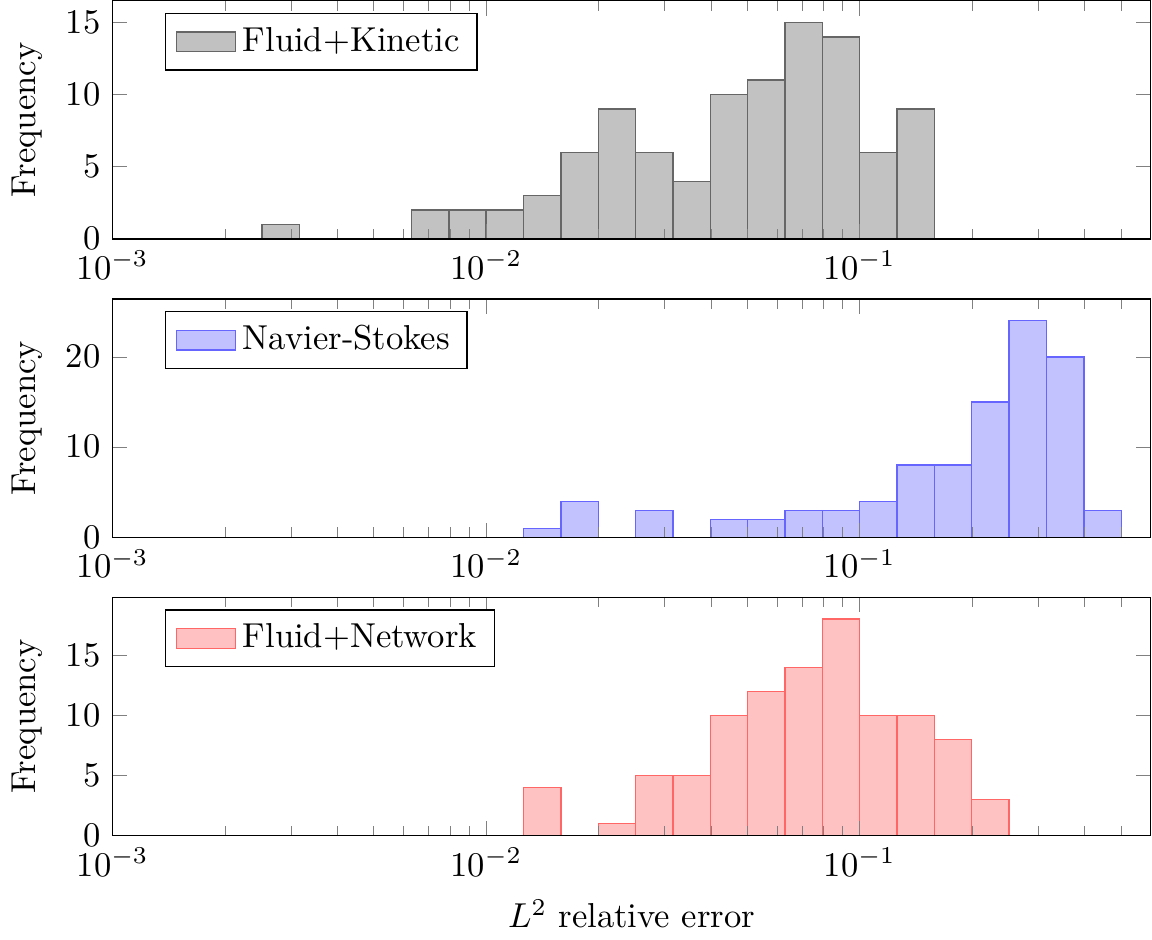}
    \caption{Distribution of the relative errors of the three fluid models on the kinetic model over 100 simulations with discontinuous initial solutions.}
    \label{fig:simu-disc}
\end{figure}

To see if that is the case, we run multiple simulations with different Knudsen numbers and different initial solutions, similar to those used to build the dataset, except that each function (density, velocity and temperature) can be multiplied by a function
$$
x \in [0, 2\pi] \mapsto \left\{ \begin{array}{cl}
\frac{c}{\pi}(x-x_d) + (1+c), & \text{if } x < x_d, \\
\frac{c}{\pi}(x-x_d) + (1-c), & \text{if } x \geqslant x_d,
\end{array}\right.
$$
adding a discontinuity at a random location $x_d\in[0,2\pi]$ with a random amplitude $c\in[-1,1]$. With these initial solutions, no simulation failed to reach $t=8$, and as can be seen Figure \ref{fig:simu-disc}, the results are comparable to those presented in Section \ref{sec:res-fluid}.

\section{Conclusion}

We construct a fluid closure for the Vlasov-Poisson dynamics based on V-net neural network and supervised learning from kinetic simulations. Slicing process of the data is introduced in order to manage meshes of different sizes. Several data processing have been also designed to improve the quality and regularity of the heat flux estimation. The numerical results show that the closure predicts the heat flux with a uniform relative error on the Knudsen interval $[0.01, 1]$, while the Navier-Stokes closure does not as expected. We also observe that the prediction is better at the beginning of the numerical simulations where the distribution function farthest from the equilibria set and so the real heat flux is larger. Surprisingly, we numerically observe that the neural network closure does not introduce instabilities when inserted in the fluid simulations, provided, however, that the outputs are regularized. Finally, the closure is quasi-optimal as the relative error between a full kinetic and a Fluid+Network closure behaves like the one between a full kinetic and a Fluid+Kinetic closure. 

This work raises several issues from the numerical point of view. The first question is its efficiency in terms of computing time or its energy cost. Presently, the Fluid+Network model requires about the same computing time as the kinetic model. However, we can hope that it will be more efficient in higher dimension, as the number of computations for a V-Net with $\ell$ levels, depth $d$ and kernel size $p$ in dimension 1 is about $O(2^\ell d^2 p N)$, while it is about $O(\ell d^2 p^2 N^2)$ in dimension 2 and $O(d^2 p^3 N^3)$ in dimension 3. Compared with the $O(N^m N_v^m)$ computations required for the kinetic model where $m$ denotes the space dimension, it increases much more slowly since $p << N_v$.
Moreover, the neural network approach can greatly benefit from GPU parallelism.

The stability of the Fluid+Network closure is also an important issue. Despite the good numerical results obtained, a mathematical guaranty is lacking. Constructing neural networks closure ensuring such stability remains to be done. 

Finally, we plan to apply this method to both the Vlasov-Poisson system in higher dimension and also to other closure problems arising in the MHD or gyrofluid design.

\bibliographystyle{plain}
\bibliography{ref}

\appendix

\section{Numerical scheme for the kinetic model}
\label{sec:vlasov-num}

In this appendix we describe the numerical method used to solve the Vlasov-Poisson equations (\ref{eq:vlasov})-(\ref{eq:poisson}) resulting from the kinetic model. This numerical method is used to produce data that can in turn be used by the neural network to interpolate the heat flux. It also serves as a reference to compute the error of the other methods, allowing us to compare them. For better readability, let us remind the one dimensional Vlasov-Poisson equations:
\begin{align*}
    \partial_t f
    +\ v\, \partial_x f
    -\ E\, \partial_v f
    =\ \frac{1}{\varepsilon} ( M(f) - f ),\\
    E = - \partial_x \phi,\quad \partial_{xx} \phi = \rho - \int_0^L \rho\,dx.
\end{align*}
The spatial domain is given by $[0,L]$, where $L > 0$ is the spatial length. For numerical purpose, the velocity domain is restricted to the bounded interval $[-v_{\max}, v_{\max} ]$. Thus we complement the equation with the following boundary conditions:
$$(\pm E)^-\, f =0,\quad \text{at } v = \pm v_{\max}.$$ 
We consider a time discretization $(t^n)_n$ with variable time step $\Delta t$ and a discretized phase space $(x_i, v_j)_{i,j}$ with constant steps $\Delta x$ and $\Delta v$ respectively. The number of discretization points in space (resp. in velocity) is denoted $N_x$ (resp. $N_v$). We denote by $f^n_{i,j}$ the approximation
\[ f^n_{i,j} \simeq f(x_i, v_j, t^n). \]
We also use the notations $\mathbf{f}^n$ for the matrix $(f^n_{i,j})_{i,j}$, $\mathbf{f}^n_i$ for the vector $(f^n_{i,j})_{j}$ and  $\mathbf{f}^n_j$ for the vector $(f^n_{i,j})_{i}$.

\subsection{Time discretization}

To solve the Vlasov-Poisson equations over the time interval $[t^n, t^n+1]$ we use a splitting between three stages: 
\begin{enumerate}
    \item Compute the electric field at time $t^n$ by solving the Poisson system: 
    $$E = - \partial_x \phi \quad, -\partial_{xx} \phi = \frac{1}{L} \int_0^L \rho dx,$$
    \item Transport the distribution function by solving the Vlasov equation over the time inverval $[t^n, t^n+1]$: 
    $$\partial_t f + v\partial_x f - E\partial_v f =0,$$
    \item Update the distribution function by taking into account the collision operator:
    $$\partial_t f = \frac{1}{\varepsilon}(M(f)-f).$$
\end{enumerate}
The first two stages rely on the space discretization and are discussed in the next section. For the third operator with a stiff source term, we use an implicit scheme:
$$
\frac{f^{n+1}-f^n}{\Delta t}=\frac{1}{\varepsilon}(M(f^{n+1})-f^{n+1}).
$$
Knowing that the fluid quantities $\rho$, $u$, $T$ are preserved by this operator, we have $M(f^{n+1}) = M(f^{n})$, so the scheme can be rewritten as:
\begin{equation}
  \label{eq:relax}
  f^{n+1} = f^n + \omega\, (M(f^n)-f^n),\quad \text{with }
  \omega = \frac{\Delta t}{\Delta t + \varepsilon}.
\end{equation}

\subsection{Spatial discretization}

For the spatial discretization, we propose a method introduced in \cite{PhamHelluy, helluy}. First we discretize in velocity with a centered finite difference scheme. After discretization, we get the following hyperbolic system:
\begin{equation}
  \label{eq:hyper}
  \frac{\mathbf{f}^{n+1}-\mathbf{f}^n}{\Delta t}+ \Lambda \partial_x \mathbf{f}^n + E B(\mathbf{f}^n) = 0
\end{equation}
with $\Lambda$ the diagonal matrix of velocities and $B(f^n)$ a vector given by
\begin{align*}
&B(f^n)_j = -\frac{f_{j+1}^n-f_{j-1}^n}{2\Delta v}, \quad j\in\{2, \dots, N_v-1\},\\
& B_1= -\frac12 \max(E,0) E f_1^n, \quad B_{N_v}= \frac12 \min(E,0) E f_{N_v}^n
\end{align*}
The choice of the boundary terms allows to obtain a dissipative hyperbolic system and to ensure $L^2$ stability. Then we discretize in space the hyperbolic system (\ref{eq:hyper}). We use a finite volume scheme with an upwind flux:
\begin{equation}
  \label{eq:transport}
  \frac{\mathbf{f}_i^{n+1}-\mathbf{f}_i^n}{\Delta t}+ \frac{\mathbf{f}_{i+\frac12}^n-\mathbf{f}_{i-\frac12}^n}{\Delta x} - E B(\mathbf{f}^n) = 0,
\end{equation}
with $\mathbf{f}_{i+\frac12}^n=\frac12\Lambda(\mathbf{f}_{i+1}^n+\mathbf{f}_{i}^n)-\frac12 \mid\!\Lambda\!\mid (\mathbf{f}_{i+1}^n-\mathbf{f}_{i}^n)$.

The Poisson equation is solved using a classical finite difference scheme:
\begin{equation}
  \label{eq:poisson-num}
  E^n_i = - \frac{\phi_{i+1}^n - \phi_{i-1}^n}{2\Delta x},\quad  -\frac{\phi_{i+1}^n -2\phi_{i}^n-\phi_{i-1}^n}{\Delta x^2} = \rho_j -\frac{1}{N_x}\sum_i^{N}\rho_i.
\end{equation}
where the density $\rho_j$ is computed as follows: $\rho_i = \sum_{j=1}^{N_v}f_{i,j}$.

The total scheme is given by (\ref{eq:poisson-num})-(\ref{eq:transport})-(\ref{eq:relax}). The time step is chosen such as to satisfy the following stability condition:
$$ \Delta t \leqslant \min \left\{\frac{\Delta x}{v_{\max}}, \frac{\Delta v}{\max_i |E_i|}\right\}.$$ 
Note that it is only first order accurate in order to avoid dispersive oscillations in the numerical solutions. 


\section{Numerical scheme for the fluid models}
\label{sec:fluid-num}

In this section we introduce the numerical methods for solving fluid models (\ref{eq:fluid}). The time discretization depends on the considered closure. For the neural network based (Fluid+Network) or the kinetic one (Fluid+Network), we use an explicit scheme. For the Navier-Stokes closure (Navier-Stokes), an implicit scheme is required to avoid a too stringent stability condition. In one case or the other, the Poisson equation is solved at the beginning of each iteration in time to compute the electric field, exactly as in Eq. \eqref{eq:poisson-num}. This section therefore focuses on the fluid equations. The following schemes are classical. We briefly present them for the sake of completeness.

\subsection{Explicit scheme for the neural network or the kinetic closure}
\label{sec:fluid-num-explicit}

Fluid equations (\ref{eq:fluid}) can be written
\begin{equation}
\partial_t \mathbf{U} + \partial_x \mathbf{F}(\mathbf{U}) = - E \mathbf{H}(\mathbf{U}),\label{systhypberbolic}
\end{equation}
with $\mathbf{U} = ( \rho, \rho u, w )$, $\mathbf{F}(\mathbf{U}) = ( \rho u, \rho u^2 + p, wu+pu+q)$, $\mathbf{H}(\mathbf{U}) = (0, \rho, \rho u)$ and where the heat flux $q$ is given by the neural network based closure (Fluid+Network) or the one obtained from full kinetic simulations (Fluid+Kinetic). We solve it on the spatial domain $[0,L]$.

When $q = 0$ (Euler closure), equation \ref{systhypberbolic} is an hyperbolic system that can be solved using a finite volume method with a local Lax-Friedrichs numerical flux and an explicit scheme in time. The hyperbolic system has characteristic speeds equal to $u, u+c, u-c$ where $c = \sqrt{3p/\rho}$ is the sound speed. When considering an additional heat flux, we propose to use the same scheme and just add a centered approximation of the heat flux in the numerical flux, as explained below.

Like for the kinetic model, we consider a time discretization $(t^n)_n$ with variable time step $\Delta t$ and a discretized phase space $(x_i)_{i}$ with constant step $\Delta x$.  We denote by $U^n_{i}$ the approximation
\[ \mathbf{U}^n_{i} \simeq \mathbf{U}(x_i, t^n). \]
The finite volume scheme results in the following formula:
\[ \frac{\mathbf{U}^{n+1}_i - \mathbf{U}^n_i}{\Delta t} + \frac{\mathbf{F}(\mathbf{U})^n_{i+\frac{1}{2}} - \mathbf{F}(\mathbf{U})^n_{i-\frac{1}{2}}}{\Delta x} = - \mathbf{H}(\mathbf{U})^n_i E^n_i,\]
with
\[ \mathbf{F}(\mathbf{U})^n_{i+\frac{1}{2}} = \frac{1}{2}( \mathbf{F}(\mathbf{U})^n_{i+1} + \mathbf{F}(\mathbf{U})^n_i) - \frac{S^n_{i+1/2}}{2}(\mathbf{U}^n_{i+1} - \mathbf{U}^n_i),\]
where
\[ S^n_{i+\frac{1}{2}} = \max(|u^n_i| + c^n_i, |u^n_{i+1}| + c^n_{i+1}) \quad\text{and}\quad c^n_i = \sqrt{\frac{3p^n_i}{\rho^n_i}}.\]
Quantities $S^n_{i+\frac{1}{2}}$ are chosen to be larger than the maximum characteristic speed of the hyperbolic system, as $u^n_i$ is the local speed of particles and $c^n_i$ the local sound speed. Comparing with the classical scheme for Euler system ($q=0$), the numerical flux for the momentum have an additional term equal to $\frac{1}{2}(q^n_{i+1} + q^n_i)$. The time step is chosen such as to satisfy the CFL stability condition:
\begin{equation}
\left(\underset{i}{\max}\ S^n_{i+\frac{1}{2}}\right) \Delta t = \frac{1}{2} \Delta x. \label{eq:CFLcond}
\end{equation}  We refer to \cite{Leveque} for more details on finite volume methods.

For the Fluid+Network method, the heat flux $q$ involved in the flux term is computed from $\varepsilon$, $\rho$, $u$ and $T$ at each iteration. For the Fluid+Kinetic method, it is obtained from an underlying kinetic simulation, using the same initial condition at $t=0$. Note that the stability condition \eqref{eq:CFLcond} does not take into account the non-zero heat flux. However, as observed in Section \ref{sec:res-fluid}, this term does not result in an additional stability condition for the Fluid+Kinetic and the Fluid+Network method, provided for the latter that the heat flux is sufficiently smoothed out. 

\subsection{Semi-implicit scheme for the Navier-Stokes closure}
\label{sec:ns-num}

With the Navier-Stokes closure $q=-\frac{3}{2}\varepsilon p\partial_x T$, the first two equations remain the same as above and are solved using the same explicit finite volume method, while the third one reads
\[ \partial_t w + \partial_x (w u + pu ) - \frac{3}{2} \varepsilon \partial_x( p \partial_x T) = - E \rho u. \]
To solve this equation we use the relation $w = \frac{1}{2}\rho u^2 + \frac{1}{2}\rho T$ to turn it into an equation on $T$, as $\rho$ and $u$ are known after solving the first two equations. We use a finite difference approximation for the term $\partial_x( p \partial_x T)$, and $T^{n+1}$ can then be computed by solving the following linear system:
\begin{align*}
    &\frac{\frac{1}{2}\rho^{n+1}_i (u^{n+1}_i)^2 + \frac{1}{2} \rho^{n+1}_i T^{n+1}_i - \frac{1}{2}\rho^{n}_i (u^{n}_i)^2 - \frac{1}{2} \rho^{n}_i T^{n}_i}{\Delta t} \\
    &+ \frac{\mathbf{F}(\mathbf{U})[2]^n_{i+\frac{1}{2}} - \mathbf{F}(\mathbf{U})[2]^n_{i-\frac{1}{2}}}{\Delta x} \\
    &- \frac{3}{2}\varepsilon \frac{ p^n_{i+1/2}T^{n+1}_{i+1} - (p^n_{i+1/2} + p^n_{i-1/2})T^{n+1}_i + p^n_{i-1/2}T^{n+1}_{i-1} }{\Delta x^2} \\
    &= - E^n_i \rho^n_i u^n_i.
\end{align*}
Here the only unknown is $T^{n+1}$, and $\mathbf{F}(\mathbf{U})[2]$ is the third coordinate of $\mathbf{F}(\mathbf{U})$. Finally, $(w)^{n+1}$ can be computed using again the relation $w = \frac{1}{2}\rho u^2 + \frac{1}{2}\rho T$. Here the time step is chosen such as to satisfy the same CFL condition \eqref{eq:CFLcond} as the implicit Navier-Stokes term does not introduce additional stability constraint.
\vspace{3cm}

Email adresses : l.bois@math.unistra.fr, emmanuel.franck@inria.fr,

laurent.navoret@math.unistra.fr, vigon@math.unistra.fr 

\end{document}